# Some applications of the Dirichlet integrals to the summation of series and the evaluation of integrals involving the Riemann zeta function

Donal F. Connon

dconnon@btopenworld.com

1 December 2012

**Abstract**

Using the Dirichlet integrals, which are employed in the theory of Fourier series, this paper develops a useful method for the summation of series and the evaluation of integrals.

## 1. Introduction

In two earlier papers [19] and [20] we considered the following identity (which is easily verified by multiplying the numerator and the denominator by the complex conjugate $(1-e^{-ix})$)

(1.1) $$\frac{1}{1-e^{ix}} = \frac{1}{2} - \frac{i}{2}\frac{\sin x}{1-\cos x} = \frac{1}{2} + \frac{i}{2}\cot(x/2)$$

Using the basic geometric series identity

$$\frac{1}{1-y} = \sum_{n=0}^{N} y^n + \frac{y^{N+1}}{1-y}$$

we obtain

(1.1a) $$\frac{1}{1-e^{ix}} = \sum_{n=0}^{N} e^{inx} + R_N(x)$$

where $$R_N(x) = \frac{e^{i(N+1)x}}{1-e^{ix}} = \frac{1}{2}e^{i(N+1)x}\left[1 + i\cot(x/2)\right] = \frac{ie^{i(N+\frac{1}{2})x}}{2\sin(x/2)}$$

Separating the real and imaginary parts of (1.1a) produces the following two identities (the first of which is called Lagrange's trigonometric identity and contains the Dirichlet kernel $D_N(x)$ which is employed in the theory of Fourier series [55, p.49])

(1.2) $$\frac{1}{2} = \sum_{n=0}^{N} \cos nx - \frac{\sin(N+1/2)x}{2\sin(x/2)}$$

(1.3) $$\frac{1}{2}\cot(x/2) = \sum_{n=0}^{N} \sin nx + \frac{\cos(N+1/2)x}{2\sin(x/2)}$$

(and these equations may obviously be generalised by substituting $\alpha x$ in place of $x$. We first multiply by $p(x)$ and, in this section, $p(x)$ is assumed to be twice continuously differentiable on $[a,b]$. Integration then gives us

$$\frac{1}{2}\int_a^b p(x)(1 + i\cot(x/2))\,dx = \sum_{n=0}^{N}\int_a^b p(x)e^{inx}\,dx + R_N$$

and we have

$$R_N = \int_a^b p(x)R_N(x)\,dx = \frac{1}{2}\int_a^b p(x)\frac{i\cos(N+1/2)x - \sin(N+1/2)x}{\sin(x/2)}\,dx$$

$$= -\frac{1}{2}\int_a^b p(x)\left\{\cos Nx + \frac{\cos(x/2)}{\sin(x/2)}\sin Nx\right\}dx$$

$$+ \frac{i}{2}\int_a^b p(x)\left\{-\sin Nx + \frac{\cos(x/2)}{\sin(x/2)}\cos Nx\right\}dx$$

Therefore, provided $\sin(x/2)$ has no zero in $[a,b]$, a weak version of the Riemann-Lebesgue lemma (which was previously employed in Eq.(2.17) in [19] where we required $p(x)$ to be twice continuously differentiable) tells us that $\lim_{N\to\infty} R_N = 0$.

This will also be the case if $\sin(a/2) = 0$, provided $p(a) = 0$. From the above we can therefore derive the following trigonometric integral identities:

(1.4) $$\frac{1}{2}\int_a^b p(x)\,dx = \sum_{n=0}^{\infty} \int_a^b p(x)\cos nx\,dx$$

(1.4a) $$\frac{1}{2}\int_a^b p(x)\cot(x/2)\,dx = \sum_{n=0}^{\infty} \int_a^b p(x)\sin nx\,dx$$

and more generally we have (by letting $x \to \alpha x$ in (1.2) and (1.3) respectively)

(1.5) $$\frac{1}{2}\int_a^b p(x)\,dx = \sum_{n=0}^{\infty} \int_a^b p(x)\cos \alpha nx\,dx$$

(1.5a) $$\frac{1}{2}\int_a^b p(x)\cot(\alpha x/2)\,dx = \sum_{n=0}^{\infty} \int_a^b p(x)\sin \alpha nx\,dx$$

Equations (1.5) and (1.5a) are valid provided (i) $\sin(\alpha x/2) \neq 0 \;\forall\; x \in [a,b]$ or, alternatively, (ii) if $\sin(\alpha a/2) = 0$ then $p(a) = 0$ also.

Similarly, using the identity

(1.6) $$\frac{1}{1+e^{ix}} = \frac{1}{2} - \frac{i}{2}\frac{\sin x}{1+\cos x} = \frac{1}{2} - \frac{i}{2}\tan(x/2)$$

we may easily prove that

(1.7) $$\frac{1}{2}\int_a^b p(x)\,dx = \sum_{n=0}^{\infty} \int_a^b p(x)(-1)^n \cos nx\,dx$$

(1.7a) $$\frac{1}{2}\int_a^b p(x)\tan(x/2)\,dx = \sum_{n=0}^{\infty} \int_a^b p(x)(-1)^{n+1} \sin nx\,dx$$

and more generally

(1.8) $$\frac{1}{2}\int_a^b p(x)\,dx = \sum_{n=0}^{\infty} \int_a^b p(x)(-1)^n \cos \alpha nx\,dx$$

(1.8a) $$\frac{1}{2}\int_a^b p(x)\tan(\alpha x/2)\,dx = \sum_{n=0}^{\infty} \int_a^b p(x)(-1)^{n+1} \sin \alpha nx\,dx$$

From (1.6) we note that the denominator is $\cos(x/2)$ and hence (1.7) and (1.7a) are only valid provided either (i) $\cos(x/2)$ has no zero in $[a,b]$ or (ii) if $\cos(a/2) = 0$, then $p(a)$ must also be zero.

Equations (1.8) and (1.8a) are valid provided (i) $\cos(\alpha x/2) \neq 0 \ \forall \ x \in [a,b]$ or, alternatively, (ii) if $\cos(\alpha a/2) = 0$ then $p(a) = 0$ also.

The following simple trigonometric identities are easily proved

(1.9) $$\cot(x/2) + \tan(x/2) = \frac{2}{\sin x}$$

(1.10) $$\cot(x/2) - \tan(x/2) = 2\cot x$$

Therefore, combining (1.5) and (1.7) produces the following identity

(1.11) $$\int_a^b \frac{p(x)}{\sin x} dx = \sum_{n=0}^{\infty} \int_a^b p(x) \sin nx \, dx - \sum_{n=0}^{\infty} \int_a^b p(x)(-1)^n \sin nx \, dx$$

which simplifies to

(1.12) $$\int_a^b \frac{p(x)}{\sin x} dx = 2 \sum_{n=0}^{\infty} \int_a^b p(x) \sin(2n+1)x \, dx$$

Similarly, using (1.10) we obtain

(1.13) $$\int_a^b p(x) \cot x \, dx = \sum_{n=0}^{\infty} \int_a^b p(x) \sin nx \, dx + \sum_{n=0}^{\infty} \int_a^b (-1)^n p(x) \sin nx \, dx$$

which simplifies to

(1.14) $$\int_a^b p(x) \cot x \, dx = 2 \sum_{n=1}^{\infty} \int_a^b p(x) \sin 2nx \, dx$$

It should be noted that in the above formulae we require either (i) both $\sin(x/2)$ and $\cos(x/2)$ have no zero in $[a,b]$ or (ii) if either $\sin(a/2)$ or $\cos(a/2)$ is equal to zero then $p(a)$ must also be zero. Condition (i) is equivalent to the requirement that $\sin x$ has no zero in $[a,b]$.

Note that (1.14) is equivalent to (1.5a) with $\alpha = 2$.

More generally we have

(1.14a) $$\int_a^b \frac{p(x)}{\sin \alpha x} dx = 2 \sum_{n=0}^{\infty} \int_a^b p(x) \sin[(2n+1)\alpha x] dx$$

(1.14b) $$\int_a^b p(x) \cot \alpha x \, dx = 2 \sum_{n=1}^{\infty} \int_a^b p(x) \sin 2\alpha n x \, dx$$

Equations (1.14a) and (1.14b) are valid provided (i) $\sin(\alpha x) \neq 0 \ \forall \ x \in [a,b]$ or, alternatively, (ii) if $\sin(\alpha a) = 0$ then $p(a) = 0$ also.

We may also generalise Poisson's integral [17, p.250] by applying the above analysis to the quotient $\frac{1}{1 \pm re^{ix}}$ and useful results may also be obtained by using $\frac{1}{1 \pm ire^{ix}}$.

## 2. Application of Dirichlet's integrals

In the previous section we required that $p(x)$ was twice continuously differentiable on the interval $[a,b]$; with the assistance of Dirichlet we now show that this condition may be significantly relaxed.

Dirichlet [26] considered the following integrals in his classical treatment of Fourier series in 1829 (when he was then 23 years old)

(2.1) $$\lim_{\mu \to \infty} \int_a^b f(t) \frac{\sin(\mu t)}{\sin t} dt = 0$$

(2.2) $$\lim_{\mu \to \infty} \int_0^b f(t) \frac{\sin(\mu t)}{\sin t} dt = \frac{\pi}{2} f(0+)$$

where $0 < a < b < \pi$. Dirichlet showed that these limits are valid if $f(x)$ is continuous on $[0,b]$ and $f(x)$ only has a finite number of maxima and minima on $[0,b]$. Jordan [37] subsequently showed in 1881 that the less restrictive condition that $f(x)$ be of bounded variation on $[0,b]$ is sufficient. Proofs may be found in [5, p.314] and [17, p.219].

We recall from [5, p.128] that if $f(x)$ is monotonic on $[a,b]$, then $f(x)$ is of bounded variation on $[a,b]$. In addition, we also note that if $f(x)$ is continuous on $[a,b]$ and if $f'(x)$ exists and is bounded on the open interval, then $f(x)$ is of bounded variation on $[a,b]$.

Even though Jordan wanted to relax the smoothness requirement for the pointwise convergence of Fourier series, it subsequently turned out that functions of bounded variation are in fact fairly smooth anyway; Lebesgue ([50, p.356] and [7, p.192]) proved that any monotonic function on a closed interval is differentiable almost everywhere, and that its derivative is integrable.

With the substitution $t = \alpha x$ in (2.1) and (2.2) we see that

$$(2.3) \qquad \lim_{\mu \to \infty} \int_{\alpha a}^{\alpha b} f(\alpha x) \frac{\sin(\mu \alpha x)}{\sin \alpha x} dx = 0$$

or equivalently

$$(2.4) \qquad \lim_{\mu \to \infty} \int_{a'}^{b'} g(x) \frac{\sin(\mu \alpha x)}{\sin \alpha x} dx = 0$$

where $0 < a' < b' < \pi$.

We may write (1.2) as

$$(2.5) \qquad \frac{1}{2} = \sum_{n=0}^{N} \cos 2\alpha n x - \frac{\sin(2N+1)\alpha x}{2 \sin \alpha x}$$

and, after multiplying this by a function $f(x)$ of bounded variation on $[a,b]$, we easily see that

$$(2.6) \qquad \frac{1}{2} \int_a^b f(x) dx = \sum_{n=0}^{N} \int_a^b f(x) \cos 2\alpha n x \, dx - \int_a^b f(x) \frac{\sin(2N+1)\alpha x}{2 \sin \alpha x} dx$$

Letting $N \to \infty$ we obtain

$$(2.7) \qquad \frac{1}{2} \int_a^b f(x) dx = \sum_{n=0}^{\infty} \int_a^b f(x) \cos 2\alpha n x \, dx$$

provided $0 < \alpha a < \alpha b < \pi$. This is a generalised version of (1.4) above and, in particular, it should be noted that we no longer require $f(x)$ to be twice continuously differentiable on $[a,b]$.

With $\alpha = \pi$ we may write (2.7) as

$$(2.8) \qquad \frac{1}{2} \int_a^b f(x) dx = \sum_{n=0}^{\infty} \int_a^b f(x) \cos 2\pi n x \, dx$$

where $0 < a < b < 1$.

We now consider the case where (2.2) applies. With the substitution $t = \alpha x$ in (2.2) we see that

(2.9) $$\lim_{\mu \to \infty} \alpha \int_0^{ab} f(\alpha x) \frac{\sin(\mu \alpha x)}{\sin \alpha x} dx = \frac{\pi}{2} f(0+)$$

or equivalently

(2.10) $$\lim_{\mu \to \infty} \int_0^{b'} g(x) \frac{\sin(\mu \alpha x)}{\sin \alpha x} dx = \frac{\pi}{2\alpha} g(0+)$$

where $0 < b' < \pi$.

Letting $N \to \infty$ in (2.6) and using (2.10) we then obtain

$$\frac{1}{2} \int_0^b f(x) dx = \sum_{n=0}^{\infty} \int_0^b f(x) \cos 2\alpha nx \, dx - \frac{\pi}{4\alpha} f(0+)$$

where $0 < \alpha b < \pi$.

With $\alpha = \pi$ we may write this as

(2.11) $$\frac{1}{2} f(0+) = \int_0^b f(x) dx + 2 \sum_{n=1}^{\infty} \int_0^b f(x) \cos 2\pi nx \, dx$$

where $0 < b < 1$.

There is clearly a connection between (2.11) and the Poisson summation formula in the form reported in Ramanujan's Notebooks [10, Part II, p.252]. If $f(x)$ is a continuous function of bounded variation on $[a,b]$, then

(2.12) $$\sum_{a \leq n \leq b}^{\#} f(n) = \int_a^b f(x) dx + 2 \sum_{n=1}^{\infty} \int_a^b f(x) \cos 2\pi nx \, dx$$

where the # on the summation sign on the left-hand side indicates that if $a$ or $b$ is an integer, then only $\frac{1}{2} f(a)$ or $\frac{1}{2} f(b)$, respectively, is counted.

The equivalence is shown below.

As noted by Bromwich [15, p.492] we consider the integral

$$\int_0^1 f(x) \frac{\sin(2N+1)\pi x}{\sin \pi x} dx = \int_0^{1/2} f(x) \frac{\sin(2N+1)\pi x}{\sin \pi x} dx + \int_{1/2}^1 f(x) \frac{\sin(2N+1)\pi x}{\sin \pi x} dx$$

and make the substitution $x = 1-t$ in the second part. This immediately gives us

$$\int_0^1 f(x)\frac{\sin(2N+1)\pi x}{\sin \pi x}dx = \int_0^{1/2}[f(x)+f(1-x)]\frac{\sin(2N+1)\pi x}{\sin \pi x}dx$$

and employing (2.10) we obtain

$$\lim_{N\to\infty}\int_0^1 f(x)\frac{\sin(2N+1)\pi x}{\sin \pi x}dx = \frac{1}{2}[f(0+)+f(1-)]$$

Therefore we have

(2.13) $\quad \frac{1}{2}[f(0+)+f(1-)] = \int_0^1 f(x)dx + 2\sum_{n=1}^\infty \int_0^1 f(x)\cos 2\pi nx\, dx$

Similarly we have

$$\int_0^2 f(x)\frac{\sin(2N+1)\pi x}{\sin \pi x}dx = \int_0^1 f(x)\frac{\sin(2N+1)\pi x}{\sin \pi x}dx + \int_1^2 f(x)\frac{\sin(2N+1)\pi x}{\sin \pi x}dx$$

and we write

$$\int_1^2 f(x)\frac{\sin(2N+1)\pi x}{\sin \pi x}dx = \int_1^{3/2} f(x)\frac{\sin(2N+1)\pi x}{\sin \pi x}dx + \int_{3/2}^2 f(x)\frac{\sin(2N+1)\pi x}{\sin \pi x}dx$$

The substitution $x = 1+t$ gives us

$$\int_1^{3/2} f(x)\frac{\sin(2N+1)\pi x}{\sin \pi x}dx = \int_0^{1/2} f(1+t)\frac{\sin(2N+1)\pi t}{\sin \pi t}dt$$

and with $x = 2-t$ we obtain

$$\int_{3/2}^2 f(x)\frac{\sin(2N+1)\pi x}{\sin \pi x}dx = \int_0^{1/2} f(2-t)\frac{\sin(2N+1)\pi t}{\sin \pi t}dt$$

Therefore, using (2.10) we see that

$$\lim_{N\to\infty}\int_0^2 f(x)\frac{\sin(2N+1)\pi x}{\sin \pi x}dx = \frac{1}{2}[f(0+)+f(1-)+f(1+)+f(2-)]$$

More generally we see that

$$\lim_{N\to\infty}\int_0^m f(x)\frac{\sin(2N+1)\pi x}{\sin \pi x}\,dx = \frac{1}{2}[f(0+)+f(m-)] + \frac{1}{2}\sum_{n=1}^{m-1}[f(n-)+f(n+)]$$

and, if $f(x)$ is continuous, this may be written as

$$\lim_{N\to\infty}\int_0^m f(x)\frac{\sin(2N+1)\pi x}{\sin \pi x}\,dx = \frac{1}{2}[f(0)+f(m)] + \sum_{n=1}^{m-1} f(n)$$

$$= -\frac{1}{2}[f(0)+f(m)] + \sum_{n=0}^{m} f(n)$$

Then using (2.6) we obtain a version of the Poisson summation formula

(2.14) $$\sum_{n=0}^{m} f(n) = \int_0^m f(x)\,dx + \frac{1}{2}[f(0)+f(m)] + 2\sum_{n=1}^{\infty}\int_0^m f(x)\cos 2n\pi x\,dx$$

Letting $m \to \infty$ and assuming that $\lim_{n\to\infty} f(n) = 0$ (which is of course required for the convergence of the infinite series) we obtain

(2.15) $$\frac{1}{2}f(0) + \sum_{n=1}^{\infty} f(n) = \int_0^{\infty} f(x)\,dx + 2\sum_{n=1}^{\infty}\int_0^{\infty} f(x)\cos 2\pi nx\,dx$$

as previously noted by Guinand [31].

Another form of the Poisson summation formula is shown below [31]

$$\sum_{n=-\infty}^{+\infty} f(n) = \sum_{n=-\infty}^{+\infty}\int_{-\infty}^{+\infty} f(x)\cos 2\pi nx\,dx$$

which may be formally obtained by letting $n \to -n$ in (2.14). An equivalent representation is

$$\sum_{n=-\infty}^{+\infty} f(n) = \sum_{n=-\infty}^{+\infty}\int_{-\infty}^{+\infty} f(x)e^{2\pi inx}\,dx$$

Various proofs exist for the Poisson formula: these include, inter alia, those given by Apostol [5, p.332], Guinand [31], Ivić [36, p.490], Mordell [42] and Wilton [53].

We now illustrate the intimate connection with the Euler-Maclaurin summation formula. Cast in its simplest form, the Euler-Maclaurin summation formula may be written as [38, p.521]

$$\sum_{n=0}^{m} f(n) = \int_0^m f(x)\,dx + \frac{1}{2}[f(0)+f(m)] + \int_0^m \left(x-[x]-\frac{1}{2}\right) f'(x)\,dx$$

and with $P_1(x) = x-[x]-\frac{1}{2}$ we have the familiar Fourier series $P_1(x) = -\sum_{n=1}^{\infty} \frac{\sin 2n\pi x}{n\pi}$.

Therefore we have

$$\int_0^m \left(x-[x]-\frac{1}{2}\right) f'(x)\,dx = -\int_0^m \sum_{n=1}^{\infty} \frac{\sin 2n\pi x}{n\pi} f'(x)\,dx$$

and with integration by parts this becomes

$$= -\sum_{n=1}^{\infty} \frac{\sin 2n\pi x}{n\pi} f(x)\bigg|_0^m + 2\int_0^m \sum_{n=1}^{\infty} f(x)\cos 2n\pi x\,dx$$

$$= 2\int_0^m \sum_{n=1}^{\infty} f(x) \cos 2n\pi x\,dx$$

Hence, assuming that interchanging the order of integration and summation is valid, we obtain the version of the Poisson summation formula in (2.15) above

$$\frac{1}{2} f(0) + \sum_{n=1}^{\infty} f(n) = \int_0^{\infty} f(x)\,dx + 2 \sum_{n=1}^{\infty} \int_0^{\infty} f(x)\cos 2\pi nx\,dx$$

□

Letting $f(x) \to f(x)\cos \pi x$ in (2.15) results in

$$\frac{1}{2} f(0) + \sum_{n=1}^{\infty} (-1)^n f(n) = \int_0^{\infty} f(x)\cos \pi x\,dx + 2\sum_{n=1}^{\infty} \int_0^{\infty} f(x) \cos 2\pi nx \cos \pi x\,dx$$

It is easily seen that

$$2\sum_{n=1}^{\infty} \int_0^{\infty} f(x)\cos 2\pi nx \cos \pi x\,dx = \sum_{n=1}^{\infty} \int_0^{\infty} f(x)\cos(2n+1)\pi x\,dx + \sum_{n=1}^{\infty} \int_0^{\infty} f(x)\cos(2n-1)\pi x\,dx$$

$$= \int_0^\infty f(x)\cos \pi x\, dx + 2\sum_{n=1}^\infty \int_0^\infty f(x)\cos(2n+1)\pi x\, dx$$

and hence we obtain

(2.16) $$\frac{1}{2}f(0) + \sum_{n=1}^\infty (-1)^n f(n) = 2\sum_{n=0}^\infty \int_0^\infty f(x)\cos(2n+1)\pi x\, dx$$

where the summation on the right-hand side now starts at $n=0$.

An alternative derivation of (2.16) is shown below. With $f(x) \to f(2x)$ we have

$$\frac{1}{2}f(0) + \sum_{n=1}^\infty f(2n) = \int_0^\infty f(2x)\, dx + 2\sum_{n=1}^\infty \int_0^\infty f(2x)\cos 2\pi nx\, dx$$

$$= \frac{1}{2}\int_0^\infty f(u)\, du + \sum_{n=1}^\infty \int_0^\infty f(u)\cos \pi nu\, du$$

Since, for suitably convergent series, $2\sum_{n=1}^\infty a_{2n} - \sum_{n=1}^\infty a_n = \sum_{n=1}^\infty (-1)^n a_n$ we have

$$\frac{1}{2}f(0) + 2\sum_{n=1}^\infty f(2n) - \sum_{n=1}^\infty f(n) = 2\sum_{n=1}^\infty \int_0^\infty f(x)[\cos \pi nx - \cos 2\pi nx]\, dx$$

$$= 2\sum_{n=0}^\infty \int_0^\infty f(x)\cos(2n+1)\pi x\, dx$$

and hence we derive (2.16) again.

### 3. Applications of the Poisson summation formula

Previously in [20] we gave many applications of the basic summation formulae stated in section 1 above. Some further applications of the more generalised versions are given below.

### 3.1 Digamma function

We consider $f(x) = \psi(a+x)$ for $\operatorname{Re}(a) > 0$ where $\psi(x)$ is the digamma function, being the logarithmic derivative of the gamma function, i.e. $\psi(x) = \dfrac{d}{dx}\log \Gamma(x)$.

We have shown in [] (this corrects the entry in [30, p.652, 6.467 2]) that

(3.1.1) $\quad \int_0^1 \psi(x+a)\cos 2n\pi x\, dx = \sin(2n\pi a)si(2n\pi a) + \cos(2n\pi a)Ci(2n\pi a)$

where $si(x)$ and $Ci(x)$ are the sine and cosine integrals defined [30, p.878] by

(3.1.2) $\quad si(x) = -\int_x^\infty \frac{\sin t}{t} dt$

and for $x > 0$

(3.1.3) $\quad Ci(x) = -\int_x^\infty \frac{\cos t}{t} dt = \gamma + \log x + \int_0^x \frac{\cos t - 1}{t} dt$

where $\gamma$ is Euler's constant.

Note that there is a slightly different sine integral $Si(x)$ which is defined in [30, p.878] and also in [1, p.231] by

(3.1.4) $\quad Si(x) = \int_0^x \frac{\sin t}{t} dt$

We have

$$si(x) = -\int_x^\infty \frac{\sin t}{t} dt = \int_0^x \frac{\sin t}{t} dt - \int_0^\infty \frac{\sin t}{t} dt$$

and using the well-known integral from Fourier series analysis [5, p.286]

(3.1.5) $\quad \dfrac{\pi}{2} = \int_0^\infty \frac{\sin t}{t} dt$

we therefore see that the two sine integrals are intimately related by

(3.1.6) $\quad si(x) = Si(x) - \dfrac{\pi}{2}$

It is easily seen that

$$\int_0^1 \psi(x+a)\, dx = \log \Gamma(1+a) - \log \Gamma(a)$$
$$= \log a$$

Since $\psi(1+x)$ is monotonic on $[0,1]$, we see that $\psi(1+x)$ is therefore of bounded variation on $[0,1]$ and hence we may apply (2.13) to obtain

$$\frac{1}{2}[\psi(a)+\psi(1+a)] = \log a + 2\sum_{n=1}^{\infty}[\cos(2n\pi a)Ci(2n\pi a)+\sin(2n\pi a)si(2n\pi a)]$$

Since $\psi(1+a) = \psi(a) + \frac{1}{a}$ this may be written as

(3.1.7) $\quad \psi(a) = \log a - \frac{1}{2a} + 2\sum_{n=1}^{\infty}[\cos(2n\pi a)Ci(2n\pi a)+\sin(2n\pi a)si(2n\pi a)]$

which appears in Nörlund's book [44, p.108]. Letting $a = 1$ results in

(3.1.8) $\quad \frac{1}{2} - \gamma = 2\sum_{n=1}^{\infty} Ci(2n\pi)$

and this corrects the corresponding formula given by Nielsen [43, p.80].

**3.2 Log gamma function**

Here we let $f(x) = \log \Gamma(a+x)$ so that (2.13) becomes

$$\frac{1}{2}[\log \Gamma(a)+\log \Gamma(1+a)] = \int_0^1 \log \Gamma(a+x)\,dx + 2\sum_{n=1}^{\infty}\int_0^1 \log \Gamma(a+x)\cos 2nx\,dx$$

We showed in [23] that for $a > 0$ and $n \geq 1$

(3.2.1) $\quad \int_0^1 \log \Gamma(x+a)\cos 2n\pi x\,dx = -\frac{1}{2n\pi}\left[-\sin(2n\pi a)Ci(2n\pi a)+\cos(2n\pi a)si(2n\pi a)\right]$

(which corrects the entry reported in [30, p.650, 6.443.3]) and hence we obtain

$$\frac{1}{2}[\log \Gamma(a)+\log \Gamma(1+a)] = \int_0^1 \log \Gamma(a+x)\,dx + \frac{1}{\pi}\sum_{n=1}^{\infty}\frac{1}{n}\left[\sin(2n\pi a)Ci(2n\pi a)-\cos(2n\pi a)si(2n\pi a)\right]$$

We designate $F(a)$ by

$$F(a) = \int_0^1 \log \Gamma(a+x)\,dx$$

and differentiation gives us

$$F'(a) = \int_0^1 \psi(a+x)\,dx = \log \Gamma(a+1) - \log \Gamma(a) = \log a$$

Integration then results in

$$F(u) - F(0) = u \log u - u$$

and using Raabe's integral [49, p.18]

$$F(0) = \int_0^1 \log \Gamma(x)\,dx = \frac{1}{2}\log(2\pi)$$

we see that

$$\int_0^1 \log \Gamma(a+x)\,dx = \frac{1}{2}\log(2\pi) + a \log a - a$$

as reported in [49, p.207].

Combining the above results in

(3.2.2)
$$\log \Gamma(a) = \frac{1}{2}\log(2\pi) + \left(a - \frac{1}{2}\right)\log a - a + \frac{1}{\pi}\sum_{n=1}^{\infty}\frac{1}{n}\left[\sin(2n\pi a)Ci(2n\pi a) - \cos(2n\pi a)si(2n\pi a)\right]$$

We note that (3.1.7) may be obtained by differentiating (3.2.2).

With $a = 1$ we obtain

(3.2.3) $$\sum_{n=1}^{\infty}\frac{si(2n\pi)}{n\pi} = \frac{1}{2}\log(2\pi) - 1$$

which was given by Nielsen [43, p.79]. With $a = 1/2$ we obtain

(3.2.4) $$\sum_{n=1}^{\infty}\frac{(-1)^n si(n\pi)}{n\pi} = \frac{1}{2}(1 - \log 2)$$

### 3.3 Logarithm function

We note that $f(x) = x \log x$ does not meet the conditions required for (1.5) on the interval $[0,1]$; however the conditions required for (2.13) are satisfied.

Integration by parts gives us

$$\int x \log x \cos ax \, dx = \frac{-Ci(ax) + ax \log x \sin ax + (1 + \log x) \cos ax}{a^2}$$

so that we have the definite integral

$$\int_0^1 x \log x \cos 2\pi nx \, dx = \frac{-Ci(2\pi n) + \lim_{x \to 0}[Ci(ax) - \cos(ax) \log x]}{(2\pi n)^2}$$

Using the definition of $Ci(x)$ in (3.3) we see that

$$Ci(ax) - \cos(ax) \log x = \gamma + \log ax - \cos(ax) \log x + \int_0^{ax} \frac{\cos t - 1}{t} dt$$

$$= \gamma + \log a - \log x [\cos(ax) - 1] + \int_0^{ax} \frac{\cos t - 1}{t} dt$$

We consider the limit

$$\lim_{x \to 0} \log x [\cos(ax) - 1] = \lim_{x \to 0} x \log x \frac{\cos(ax) - 1}{x}$$

$$= \lim_{x \to 0} x \log x \cdot \lim_{x \to 0} \frac{\cos(ax) - 1}{x}$$

and using L'Hôpital's rule we obtain

$$\lim_{x \to 0} \frac{\cos(ax) - 1}{x} = 0$$

which shows that

$$\lim_{x \to 0} [\cos(ax) - 1] \log x = 0$$

Then taking the limit as $x \to 0$ we obtain

$$\lim_{x \to 0} [Ci(ax) - \cos(ax) \log x] = \gamma + \log a$$

Hence we have

$$\int_0^1 x \log x \cos 2\pi nx\, dx = \frac{-Ci(2\pi n) + \gamma + \log(2\pi n)}{(2\pi n)^2}$$

and the integral $\int_0^1 x \log x\, dx = -1/4$ is elementary.

Therefore, with $f(x) = x \log x$ in (2.13) we obtain

$$\frac{1}{4} = 2 \sum_{n=1}^{\infty} \frac{-Ci(2\pi n) + \gamma + \log(2\pi n)}{(2\pi n)^2}$$

so that

(3.3.1) $$\sum_{n=1}^{\infty} \frac{Ci(2n\pi)}{n^2} = \varsigma(2)[\gamma + \log(2\pi)] - \varsigma'(2) - \frac{1}{2}\pi^2$$

**3.4 An integral due to Ramanujan**

With $f(u) = \psi(1+u) - \log(1+u) + \frac{1}{2(1+u)}$ in (2.15) we have

(3.4.1)
$$\frac{1}{2}\left[\frac{1}{2} - \gamma\right] + \sum_{n=1}^{\infty}\left[\psi(1+n) - \log(1+n) + \frac{1}{2(1+n)}\right] = \int_0^{\infty}\left[\psi(1+u) - \log(1+u) + \frac{1}{2(1+u)}\right] du$$

$$+ 2\sum_{n=1}^{\infty} \int_0^{\infty}\left[\psi(1+u) - \log(1+u) + \frac{1}{2(1+u)}\right] \cos 2\pi nu\, du$$

The left-hand side of (3.4.1) may be written as

$$\frac{1}{2}\left[\frac{1}{2} - \gamma\right] + \sum_{n=1}^{\infty}\left[\psi(1+n) - \log(1+n) + \frac{1}{2(1+n)}\right] = \sum_{n=1}^{\infty}\left[\psi(n) - \log n + \frac{1}{2n}\right] - \frac{1}{2}\left[\frac{1}{2} - \gamma\right]$$

and we split the other integral on the right-hand side into three components as shown below

$$\int_0^{\infty}\left[\psi(1+u) - \log(1+u) + \frac{1}{2(1+u)}\right] \cos 2\pi nu\, du$$

$$= \int_0^{\infty}[\psi(1+u) - \log u] \cos 2\pi nu\, du + \int_0^{\infty}[\log u - \log(1+u)] \cos 2\pi nu\, du$$

$$+\frac{1}{2}\int_0^\infty \frac{\cos 2\pi nu}{1+u}du$$

Ramanujan [12] obtained the following integral

(3.4.2) $$\int_0^\infty [\psi(1+u)-\log u]\cos 2n\pi u\, du = \frac{1}{2}[\psi(1+n)-\log n]$$

which was rediscovered by Guinand [33] in 1947.

As regards the next component, integration by parts readily gives us

$$a\int [\log u - \log(1+u)]\cos au\, du = -\sin a\, Ci[a(1+u)]+\cos a\, Si[a(1+u)] - Si(au)$$

$$+[\log u - \log(1+u)]\sin au$$

With $a = 2\pi n$ this simplifies to

$$2\pi n\int [\log u - \log(1+u)]\cos 2\pi nu\, du = Si[2\pi n(1+u)] - Si(2\pi nu) + [\log u - \log(1+u)]\sin 2\pi nu$$

and we obtain the definite integral

$$\int_0^\infty [\log u - \log(1+u)]\cos 2\pi nu\, du = \frac{Si(0)-Si(2\pi n)}{2\pi n}$$

$$= -\frac{Si(2\pi n)}{2\pi n} = -\frac{si(2\pi n)}{2\pi n} - \frac{1}{4n}$$

The substitution $v = 1+u$ results in

$$\int \frac{\cos au}{1+u}du = \cos a\, Ci[a(1+u)] + \sin a\, Si[a(1+u)]$$

where $Si(u) = \frac{\pi}{2} + si(u)$. Therefore we have the definite integral

$$\int_0^\infty \frac{\cos au}{1+u}du = \sin a\, Si(\infty) - [\cos a\, Ci(a) + \sin a\, Si(a)]$$

and with $a = 2\pi n$ this becomes

$$\int_0^\infty \frac{\cos 2n\pi a u}{1+u} du = -Ci(2n\pi)$$

With $I = \int_0^\infty \left[ \psi(1+u) - \log(1+u) + \frac{1}{2(1+u)} \right] du$ we then have

$$\sum_{n=1}^\infty \left[ \psi(n) - \log n + \frac{1}{2n} \right] - \frac{1}{2}\left[\frac{1}{2} - \gamma\right] = I + 2\sum_{n=1}^\infty \left[ \frac{1}{2}[\psi(1+n) - \log n] - \frac{Si(2\pi n)}{2\pi n} - \frac{1}{2}Ci(2n\pi) \right]$$

$$= I + \sum_{n=1}^\infty \left[ [\psi(1+n) - \log n] - \frac{1}{2n} - \frac{si(2\pi n)}{\pi n} - Ci(2n\pi) \right]$$

$$= I + \sum_{n=1}^\infty \left[ \psi(n) - \log n + \frac{1}{2n} \right] - \left[ \frac{1}{2}\log(2\pi) - 1 \right] - \frac{1}{2}\left[\frac{1}{2} - \gamma\right]$$

where we have employed the series

$$\sum_{n=1}^\infty Ci(2\pi n) = \frac{1}{2}\left[\frac{1}{2} - \gamma\right]$$

$$\sum_{n=1}^\infty \frac{si(2\pi n)}{\pi n} = \frac{1}{2}\log(2\pi) - 1$$

which were derived in (3.1.8) and (3.2.1) respectively.

We therefore obtain

(3.4.3) $$\int_0^\infty \left[ \psi(1+u) - \log(1+u) + \frac{1}{2(1+u)} \right] du = \frac{1}{2}\log(2\pi) - 1$$

□

We now give an alternative derivation of the integral (3.4.3). We recall Binet's second formula for $\log \Gamma(u)$ (which is derived in [51, p.250] using the Abel-Plana summation formula)

(3.4.4) $$\log \Gamma(u) = \left( u - \frac{1}{2} \right) \log u - u + \frac{1}{2}\log(2\pi) + 2\int_0^\infty \frac{\tan^{-1}(x/u)}{e^{2\pi x} - 1} dx$$

Another proof of this is also reported in [22]. This formula was also derived by Ramanujan [10, Part II, p.221] in the case where $u$ is a positive integer.

Differentiation of Binet's formula results in [51, p.251]

$$(3.4.5) \qquad \log u - \psi(u) = \frac{1}{2u} + 2\int_0^\infty \frac{x}{(u^2 + x^2)(e^{2\pi x} - 1)} dx$$

and thus we have

$$(3.4.6) \qquad \psi(t+u) - \log(t+u) + \frac{1}{2(t+u)} = -2\int_0^\infty \frac{x}{[(t+u)^2 + x^2](e^{2\pi x} - 1)} dx$$

Integration with respect to $u$ gives us

$$I(t) = \int_0^\infty \left[ \psi(t+u) - \log(t+u) + \frac{1}{2(t+u)} \right] du = -2\int_0^\infty \int_0^\infty \frac{x}{[(t+u)^2 + x^2](e^{2\pi x} - 1)} dx\, du$$

$$= -2 \int_0^\infty \frac{x\, dx}{e^{2\pi x} - 1} \int_0^\infty \frac{du}{(t+u)^2 + x^2}$$

We see that

$$\int_0^\infty \frac{du}{(t+u)^2 + x^2} = \int_t^\infty \frac{dy}{y^2 + x^2}$$

$$= \frac{1}{x} \tan^{-1}\left(\frac{y}{x}\right) \Big|_t^\infty$$

$$= \frac{1}{x}\left[\frac{\pi}{2} - \tan^{-1}\left(\frac{t}{x}\right)\right]$$

Since

$$\tan^{-1} a + \tan^{-1} b = \tan^{-1}\left(\frac{a+b}{1-ab}\right)$$

we see that

$$\tan^{-1}(x/u) = \frac{\pi}{2} - \tan^{-1}(u/x) \text{ for } u/x > 0$$

and the integral may thus be written as

$$\int_0^\infty \frac{du}{(t+u)^2 + x^2} = \frac{1}{x} \tan^{-1}\left(\frac{x}{t}\right)$$

We therefore have

$$I(t) = -2\int_0^\infty \frac{\tan^{-1}(x/t)}{e^{2\pi x}-1}dx$$

We then obtain from (3.4.4)

$$I(t) = \left(t-\frac{1}{2}\right)\log t - t + \frac{1}{2}\log(2\pi) - \log\Gamma(t)$$

so that

(3.4.7) $$\int_0^\infty\left[\psi(t+u) - \log(t+u) + \frac{1}{2(t+u)}\right]du = \left(t-\frac{1}{2}\right)\log t - t + \frac{1}{2}\log(2\pi) - \log\Gamma(t)$$

Letting $t=1$ gives us (3.4.3)

(3.4.8) $$\int_0^\infty\left[\psi(1+u) - \log(1+u) + \frac{1}{2(1+u)}\right]du = \frac{1}{2}\log(2\pi) - 1$$

It is easy to determine that

$$\int_0^N\left[\log u - \log(1+u) + \frac{1}{1+u}\right]du = N\log N - (1+N)\log(1+N) + \log(1+N)$$

$$= -N\log\left(1+\frac{1}{N}\right)$$

and, since $\log\left(1+\frac{1}{N}\right) = \frac{1}{N} + O\left(\frac{1}{N^2}\right)$, we then see that

$$\lim_{N\to\infty}\int_0^N\left[\log u - \log(1+u) + \frac{1}{1+u}\right]du = -1$$

We note that

$$\int_0^\infty\left[\psi(1+u) - \log(1+u) + \frac{1}{2(1+u)}\right]du$$

$$= \int_0^\infty\left[\psi(1+u) - \log u - \frac{1}{2(1+u)} + \log u - \log(1+u) + \frac{1}{1+u}\right]du$$

and thus we have

$$\int_0^\infty \left[\psi(1+u) - \log(1+u) + \frac{1}{2(1+u)}\right] du = -1 + \int_0^\infty \left[\psi(1+u) - \log u - \frac{1}{2(1+u)}\right] du$$

Therefore using (3.4.8) we obtain

(3.4.9) $$\int_0^\infty \left[\psi(1+u) - \log u - \frac{1}{2(1+u)}\right] du = \frac{1}{2}\log(2\pi)$$

as previously obtained by Berndt and Dixit [12]. This integral may also be evaluated in the following manner:

We see that

$$\int_0^\infty \left[\psi(1+u) - \log u - \frac{1}{2(u+1)}\right] du = \left. \log \Gamma(1+u) - u \log u + u - \frac{1}{2}\log(1+u) \right|_0^\infty$$

$$= \lim_{u \to \infty} \left[ \log \Gamma(1+u) - u \log u + u - \frac{1}{2}\log(1+u) \right]$$

$$= \lim_{u \to \infty} \left[ \log \Gamma(1+u) - u \log u + u - \frac{1}{2}\log u - \frac{1}{2}\log\left(1+\frac{1}{u}\right) \right]$$

$$= \lim_{u \to \infty} \left[ \log \Gamma(1+u) - \left(u + \frac{1}{2}\right) u \log u + u \right]$$

The integral (3.4.9) may then be obtained upon using the asymptotic expression [49, p.8]

$$\log \Gamma(1+u) \approx \left(u + \frac{1}{2}\right) \log u - u + \frac{1}{2}\log(2\pi) + \frac{1}{12u}$$

A generalised version of the integral (3.4.2) appears in [30, p.652, 6.471.3]

(3.4.10) $$\int_0^\infty [\psi(1+u) - \log u] \cos 2\pi ut \, du = \frac{1}{2}[\psi(1+t) - \log t]$$

and integration results in

$$(3.4.11) \quad \frac{1}{2\pi}\int_0^\infty \frac{\psi(1+u)-\log u}{u}\sin 2\pi uv\, du = \frac{1}{2}\left[\log\Gamma(1+v)-v\log v+v\right]$$

□

Merkle and Merkle [40] have recently shown that

$$(3.4.12)\quad \log\Gamma(x) = \frac{1}{2}[\log(2\pi)-1]-\frac{1}{2}\gamma(2x-1)+\sum_{n=0}^\infty\left[\psi(1+n)-\log(x+n)+\frac{2x-1}{2(1+n)}\right]$$

and commencing the summation at $n=1$ gives us

(3.4.13)
$$\log\Gamma(1+x) = \frac{1}{2}[\log(2\pi)-1]-\gamma\left(x+\frac{1}{2}\right)+\frac{1}{2}(2x-1)+\sum_{n=1}^\infty\left[\psi(1+n)-\log(x+n)+\frac{2x-1}{2(1+n)}\right]$$

Differentiation gives us

$$\psi(1+x) = -\gamma+1+\sum_{n=1}^\infty\left[-\frac{1}{x+n}+\frac{1}{1+n}\right]$$

$$= -\gamma+1+\sum_{n=1}^\infty\left[-\frac{1}{x+n}+\frac{1}{n}+\frac{1}{1+n}-\frac{1}{n}\right]$$

$$= -\gamma+1+\sum_{n=1}^\infty\left[-\frac{1}{x+n}+\frac{1}{n}\right]+\sum_{n=1}^\infty\left[\frac{1}{1+n}-\frac{1}{n}\right]$$

Hence we obtain the well-known formula for the digamma function [49, p.14]

$$(3.4.14)\quad \psi(1+x) = -\gamma-\sum_{n=1}^\infty\left[\frac{1}{x+n}-\frac{1}{n}\right]$$

Letting $x=1$ in (3.4.13) gives us

$$0 = \frac{1}{2}[\log(2\pi)-1]-\frac{3}{2}\gamma+\frac{1}{2}+\sum_{n=1}^\infty\left[\psi(1+n)-\log(1+n)+\frac{1}{2(1+n)}\right]$$

which may be written as

$$0 = \frac{1}{2}\log(2\pi)-\frac{3}{2}\gamma+\sum_{n=1}^\infty\left[\psi(n)-\log n+\frac{1}{2n}\right]-\psi(1)-\frac{1}{2}$$

Hence we have

(3.4.15) $$\sum_{n=1}^{\infty}\left[\psi(n)-\log n+\frac{1}{2n}\right]=\frac{1}{2}[1+\gamma-\log(2\pi)]$$

and we note that this series appears in (3.4.1).

With $x=0$ in (3.4.13) we get

(3.4.15.1) $$\sum_{n=1}^{\infty}\left[\psi(1+n)-\log n-\frac{1}{2(1+n)}\right]=1+\frac{1}{2}\gamma-\frac{1}{2}\log(2\pi)$$

The left-hand side of (3.4.15.1) may be written as

(3.4.15.2) $$\sum_{n=1}^{\infty}\left[\psi(n)-\log n+\frac{1}{2n}+\frac{1}{2n}-\frac{1}{2(1+n)}\right]=1+\frac{1}{2}\gamma-\frac{1}{2}\log(2\pi)$$

and we see that this is consistent with (3.4.15) because subtraction results in the telescoping series

$$\sum_{n=1}^{\infty}\left[\frac{1}{n}-\frac{1}{1+n}\right]=1$$

Letting $x=1/2$ in (3.4.13) gives us

(3.4.16) $$\frac{1}{2}[1-\log 2]=\sum_{n=0}^{\infty}[\psi(1+n)+\log 2-\log(2n+1)]$$

The Weierstrass expression for the gamma function may be written as [49, p.1]

(3.4.17) $$\log\Gamma(x+a)=-\log(x+a)-\gamma(x+a)+\sum_{n=1}^{\infty}\left[\log n-\log(n+a+x)+\frac{x+a}{n}\right]$$

so that

(3.4.18) $$\log\Gamma(1+x)=-\log(1+x)-\gamma(1+x)+\sum_{n=1}^{\infty}\left[\log n-\log(n+1+x)+\frac{1+x}{n}\right]$$

Subtracting (3.4.18) from (3.4.13) results in

(3.4.19) $$\sum_{n=1}^{\infty}\left[\psi(1+n)-\log(x+n)-\log n+\log(n+1+x)+\frac{2x-1}{2(1+n)}-\frac{1+x}{n}\right]$$

$$=1-\frac{1}{2}\log(2\pi)-\frac{1}{2}\gamma-x-\log(1+x)$$

and with $x = 0$ we obtain

$$\sum_{n=1}^{\infty}\left[\psi(1+n)-2\log n+\log(n+1)-\frac{1}{2(1+n)}-\frac{1}{n}\right]=1-\frac{1}{2}\gamma-\frac{1}{2}\log(2\pi)$$

which may be expressed as

$$\sum_{n=1}^{\infty}\left[\psi(1+n)-\log n-\frac{1}{2(1+n)}\right]-\sum_{n=1}^{\infty}\left[\frac{1}{n}-\log\left(1+\frac{1}{n}\right)\right]=1-\frac{1}{2}\gamma-\frac{1}{2}\log(2\pi)$$

and this concurs with the result obtained by letting $x = 0$ in (3.4.13).

$\square$

Merkle and Merkle [40] also showed that

(3.4.20) $$\sum_{n=0}^{\infty}\left[\psi(x+n)-\psi(1+n)-\frac{x-1}{1+n}\right]=(1-x)[\psi(x)+\gamma-1]$$

where a misprint in their formula has been corrected. Starting the summation at $n = 1$ gives us

(3.4.21) $$\sum_{n=1}^{\infty}\left[\psi(x+n)-\psi(1+n)-\frac{x-1}{1+n}\right]=-x[\psi(x)+\gamma]+2(x-1)$$

Differentiation results in

(3.4.22) $$\sum_{n=1}^{\infty}\left[\psi'(x+n)-\frac{1}{1+n}\right]=-x\psi'(x)-\psi(x)-\gamma+2$$

so that

$$\sum_{n=1}^{\infty}\left[\psi'(1+n)-\frac{1}{1+n}\right]=-\psi'(1)+2$$

or equivalently

(3.4.23) $$\sum_{n=1}^{\infty}\left[\psi'(n)-\frac{1}{n}\right]=1$$

which may also be expressed as

$$\sum_{n=1}^{\infty}\left[\varsigma(2,n)-\frac{1}{n}\right]=1$$

We note that (3.4.23) is consistent with the asymptotic formula found by Barnes in 1899 [49, p.23]

$$\sum_{n=1}^{m}\psi'(n)=\log m+1+\gamma+O(m^{-1})\qquad m\to\infty$$

Adding (3.4.15) and (3.4.23) gives us

(3.4.24) $$\sum_{n=1}^{\infty}\left[\psi(n)-\log n+\frac{1}{2}\psi'(n)\right]=1+\frac{1}{2}\gamma-\frac{1}{2}\log(2\pi)$$

which, as stated by Srivastava and Choi [49, p.29], may also be derived from the following representation of the Barnes double gamma function

$$G(1+x)=(2\pi)^{\frac{1}{2}x}\exp\left[-\frac{1}{2}x(x+1)-\frac{1}{2}\gamma x^2\right]\prod_{n=1}^{\infty}\frac{\Gamma(n)}{\Gamma(x+n)}\exp\left[x\psi(n)+\frac{1}{2}x^2\psi'(n)\right]$$

by taking logarithms of both sides and letting $x=1$ in the resulting equation.

□

Differentiating (3.4.22) gives us

$$\sum_{n=1}^{\infty}\psi''(x+n)=-x\psi''(x)-2\psi'(x)$$

and we see that

(3.4.25) $$\sum_{n=1}^{\infty}\psi^{(p)}(x+n)=-x\psi^{(p)}(x)-p\psi^{(p-1)}(x)$$

which was obtained in a different manner by Merkle and Merkle [40].

We need to reconcile this with the result previously found by Adamchik and Srivastava ([2] and [49, p.156]) which is valid for $|z|<1$

(3.4.26) $$\sum_{n=0}^{\infty}\psi^{(p)}(x+n)z^n=\frac{(-1)^{p+1}p!}{x^{p+1}}+\frac{\psi^{(p)}(x+1)+(-1)^p p!z^2\Phi(z,p+1,x+1)}{1-z}$$

where $\Phi(z,p,x)$ is the Hurwitz-Lerch zeta function defined by

$$\Phi(z, p, x) = \sum_{n=0}^{\infty} \frac{z^n}{(x+n)^p}$$

It should be noted that (3.4.26) has been corrected to show that the summation starts at $n = 0$ (this is because the summation in the second equation on page 157 of [49] should also have been started at $k = 0$).

We may express (3.4.26) as

$$\sum_{n=0}^{\infty} \psi^{(p)}(x+n) z^n = \frac{(-1)^{p+1} p!}{x^{p+1}} + \frac{(-1)^{p+1} p! \left[ \varsigma(p+1, x+1) - z^2 \Phi(z, p+1, x+1) \right]}{1-z}$$

and, since $\lim_{z \to 1} \left[ \varsigma(p+1, x+1) - z^2 \Phi(z, p+1, x+1) \right] = 0$, we may apply L'Hôpital's rule to obtain

$$\sum_{n=0}^{\infty} \psi^{(p)}(x+n) = \frac{(-1)^{p+1} p!}{x^{p+1}} + (-1)^{p+1} p! \lim_{z \to 1} \left[ 2z \Phi(z, p+1, x+1) + z^2 \Phi'(z, p+1, x+1) \right]$$

$$= \frac{(-1)^{p+1} p!}{x^{p+1}} + (-1)^{p+1} p! \lim_{z \to 1} \left[ 2\varsigma(p+1, x+1) + \sum_{n=0}^{\infty} \frac{n z^{n+1}}{(x+1+n)^{p+1}} \right]$$

$$= \frac{(-1)^{p+1} p!}{x^{p+1}} + (-1)^{p+1} p! \sum_{n=0}^{\infty} \frac{2+n}{(x+1+n)^{p+1}}$$

Since

$$\sum_{n=0}^{\infty} \frac{2+n}{(x+1+n)^{p+1}} = \sum_{n=0}^{\infty} \frac{x+1+n}{(x+1+n)^{p+1}} + \sum_{n=0}^{\infty} \frac{1-x}{(x+1+n)^{p+1}}$$

we see that

$$\sum_{n=0}^{\infty} \psi^{(p)}(x+n) = \frac{(-1)^{p+1} p!}{x^{p+1}} + (-1)^{p+1} p! [\varsigma(p, x+1) + (1-x)\varsigma(p+1, x+1)]$$

$$= \frac{(-1)^{p+1} p!}{x^{p+1}} - p\psi^{(p-1)}(x+1) + (1-x)\psi^{(p)}(x+1)$$

Since $\psi(1+x) = \psi(x) + \frac{1}{x}$ this may be written as

$$\sum_{n=0}^{\infty} \psi^{(p)}(x+n) = (1-x)\psi^{(p)}(x) - p\psi^{(p-1)}(x)$$

which concurs with (3.4.25).

$\square$

Equation (3.4.13) may be written as

$$\log \Gamma(1+x) = \frac{1}{2}[\log(2\pi)-1] + \frac{1}{2}(1-\gamma)(2x-1) - \gamma + \sum_{n=1}^{\infty}\left[\psi(1+n) - \log(x+n) + \frac{2x-1}{2(1+n)}\right]$$

We multiply this by $\sin 2k\pi x$ and integrate this to obtain

$$\int_0^1 \log \Gamma(1+x) \sin 2k\pi x\, dx = -\frac{1}{2k\pi}(1-\gamma) - \sum_{n=1}^{\infty}\left[\int_0^1 \log(x+n)\sin 2k\pi x\, dx + \frac{1}{2k\pi(1+n)}\right]$$

Using

$$2k\pi \int_0^1 \log(x+n)\sin 2k\pi x\, dx = Ci[2k\pi(n+1)] - Ci[2k\pi n] - \log(n+1) + \log n$$

we have the finite sum

$$2k\pi \sum_{n=1}^{N} \int_0^1 \log(x+n) \sin 2k\pi x\, dx = \sum_{n=1}^{N}\left[Ci[2k\pi(n+1)] - Ci[2k\pi n] - \log(n+1) + \log n\right]$$

and this telescopes to

$$= Ci[2k\pi(N+1)] - Ci[2k\pi] - \log(N+1)$$

We then have

$$2k\pi \sum_{n=1}^{N}\left[\int_0^1 \log(x+n)\sin 2k\pi x\, dx + \frac{1}{1+n}\right] = Ci[2k\pi(N+1)] - Ci[2k\pi] - \log(N+1) + H_{N+2} - 1$$

and thus we have the limit as $N \to \infty$

$$2k\pi \sum_{n=1}^{\infty}\left[\int_0^1 \log(x+n)\sin 2k\pi x\, dx + \frac{1}{1+n}\right] = -Ci(2k\pi) + \gamma - 1$$

Hence we obtain the known integral

$$\int_0^1 \log \Gamma(1+x)\sin 2k\pi x\, dx = \frac{Ci(2k\pi)}{2k\pi}$$

□

Merkle and Merkle [40] also showed that

(3.4.27) $\log G(x) = \frac{1}{2}[\log(2\pi)-1](x-1) - \frac{1}{2}\gamma(x-1)^2$

$$-\sum_{n=0}^{\infty}\left[\log\Gamma(x+n) - \log\Gamma(1+n) - \psi(1+n)(x-1) - \frac{(x-1)^2}{2(1+n)}\right]$$

so that

(3.4.28) $\log G(1+x) = \frac{1}{2}[\log(2\pi)-1]x - \frac{1}{2}\gamma x^2$

$$-\sum_{n=0}^{\infty}\left[\log\Gamma(1+x+n) - \log\Gamma(1+n) - \psi(1+n)x - \frac{x^2}{2(1+n)}\right]$$

It is well known that the Barnes double gamma function may be represented by [49, p.25]

(3.4.29) $\log G(1+x) = \frac{1}{2}x\log(2\pi) - \frac{1}{2}x(1+x) - \frac{1}{2}\gamma x^2 + \sum_{n=1}^{\infty}\left[\frac{1}{2n}x^2 - x + n\log\left(1+\frac{x}{n}\right)\right]$

and equating (3.4.28) with (3.4.29) gives us

$$-\frac{1}{2}x^2 + \sum_{n=1}^{\infty}\left[\frac{1}{2n}x^2 - x + n\log\left(1+\frac{x}{n}\right)\right] = -\sum_{n=0}^{\infty}\left[\log\Gamma(1+x+n) - \log\Gamma(1+n) - \psi(1+n)x - \frac{x^2}{2(1+n)}\right]$$

which may be written as

(3.4.30) $\frac{1}{2}x^2 = \sum_{n=1}^{\infty}\left[\log\Gamma(x+n) - \log\Gamma(n) + n\log\left(1+\frac{x}{n}\right) - \psi(n)x - x\right]$

With $x = 1$ we get

(3.4.31) $\frac{1}{2} = \sum_{n=1}^{\infty}\left[\log n + n\log\left(1+\frac{1}{n}\right) - \psi(n) - 1\right]$

Combining this with (3.4.15) results in

$$(3.4.32) \qquad 1+\frac{1}{2}\gamma-\frac{1}{2}\log(2\pi)=\sum_{n=1}^{\infty}\left[\log\left(1+\frac{1}{n}\right)^{n}+\frac{1}{2n}-1\right]$$

Since $\lim_{n\to\infty}\log\left(1+\frac{1}{n}\right)^{n}=\log e=1$, we note that $\lim_{n\to\infty}\left[\log\left(1+\frac{1}{n}\right)^{n}+\frac{1}{2n}-1\right]=0$ and hence, as expected, the $n$th term of the series (3.4.32) approaches zero as $n\to\infty$.

We showed in [21] that

$$1-\frac{1}{2}\log(2\pi)=\sum_{n=1}^{\infty}\left[\frac{1}{2}\log\left(1+\frac{1}{n}\right)-1+\log\left(1+\frac{1}{n}\right)^{n}\right]$$

and subtraction of (3.4.32) gives us the well known result [47]

$$\gamma=\sum_{n=1}^{\infty}\left[\frac{1}{n}-\log\left(1+\frac{1}{n}\right)\right]$$

□

Letting $x=1/2$ in (3.4.21) gives us

$$\sum_{n=1}^{\infty}\left[\psi\left(n+\frac{1}{2}\right)-\psi(1+n)+\frac{1}{2(1+n)}\right]=-\frac{1}{2}\left[\psi\left(\frac{1}{2}\right)+\gamma\right]-1$$

$$=\log 2-1$$

It is easily seen that

$$\sum_{n=1}^{\infty}\left[\psi\left(n+\frac{1}{2}\right)-\log n\right]=\sum_{n=1}^{\infty}\left[\psi\left(n+\frac{1}{2}\right)-\psi(1+n)+\frac{1}{2(1+n)}\right]+\sum_{n=1}^{\infty}\left[\psi(1+n)-\log n-\frac{1}{2(1+n)}\right]$$

and substituting (3.4.15.1) we get

$$(3.4.33) \qquad \sum_{n=1}^{\infty}\left[\psi\left(n+\frac{1}{2}\right)-\log n\right]=\frac{1}{2}\gamma+\frac{1}{2}\log\frac{2}{\pi}$$

This was obtained in a different way in [25].

□

Another derivation of (3.4.33) is shown below. Differentiating (3.4.30) gives us

(3.4.34) $$x = \sum_{n=1}^{\infty}\left[\psi(x+n)-\psi(n)-\frac{x}{n+x}\right]$$

and with $x=1/2$ we have

$$\frac{1}{2} = \sum_{n=1}^{\infty}\left[\psi\left(n+\frac{1}{2}\right)-\psi(n)-\frac{1}{2n+1}\right]$$

$$= \sum_{n=1}^{\infty}\left[\psi\left(n+\frac{1}{2}\right)-\log n+\log n-\psi(n)-\frac{1}{2n+1}\right]$$

$$= \sum_{n=1}^{\infty}\left[\psi\left(n+\frac{1}{2}\right)-\log n\right]-\sum_{n=1}^{\infty}\left[\psi(n)-\log n+\frac{1}{2n+1}\right]$$

$$= \sum_{n=1}^{\infty}\left[\psi\left(n+\frac{1}{2}\right)-\log n\right]-\sum_{n=1}^{\infty}\left[\psi(n+1)-\log n-\frac{1}{n}+\frac{1}{2n+1}\right]$$

$$= \sum_{n=1}^{\infty}\left[\psi\left(n+\frac{1}{2}\right)-\log n\right]-\sum_{n=1}^{\infty}\left[\psi(n+1)-\log n-\frac{1}{2(1+n)}+\frac{1}{2(1+n)}-\frac{1}{n}+\frac{1}{2n+1}\right]$$

$$= \sum_{n=1}^{\infty}\left[\psi\left(n+\frac{1}{2}\right)-\log n\right]-\sum_{n=1}^{\infty}\left[\psi(n+1)-\log n-\frac{1}{2(1+n)}\right]$$

$$-\sum_{n=1}^{\infty}\left[\frac{1}{2(1+n)}-\frac{1}{n}+\frac{1}{2n+1}\right]$$

Using (3.4.15.1) we obtain

$$\frac{1}{2} = \sum_{n=1}^{\infty}\left[\psi\left(n+\frac{1}{2}\right)-\log n\right]-1-\frac{1}{2}\gamma+\frac{1}{2}\log(2\pi)-\sum_{n=1}^{\infty}\left[\frac{1}{2(1+n)}-\frac{1}{n}+\frac{1}{2n+1}\right]$$

We see that

$$\sum_{n=1}^{\infty}\left[\frac{1}{2(1+n)}-\frac{1}{n}+\frac{1}{2n+1}\right] = \sum_{n=1}^{\infty}\left[\frac{1}{2n+2}-\frac{1}{2n}\right]+\sum_{n=1}^{\infty}\left[\frac{1}{2n+1}-\frac{1}{2n}\right]$$

$$= \frac{1}{2}\sum_{n=1}^{\infty}\left[\frac{1}{n+1}-\frac{1}{n}\right]+\frac{1}{2}\sum_{n=1}^{\infty}\left[\frac{1}{n+\frac{1}{2}}-\frac{1}{n}\right]$$

and using (3.4.14)

$$\psi(1+x) = -\gamma - \sum_{n=1}^{\infty}\left[\frac{1}{x+n} - \frac{1}{n}\right]$$

we obtain

$$\sum_{n=1}^{\infty}\left[\frac{1}{2(1+n)} - \frac{1}{n} + \frac{1}{2n+1}\right] = -\frac{1}{2}[\psi(2)+\gamma] - \frac{1}{2}[\psi(3/2)+\gamma]$$

$$= -\frac{1}{2} + \log 2 - 1$$

and then deduce (3.4.33).

□

In passing, we note that Guinand [33] has shown that for $|\arg z| < \pi$

$$(3.4.35) \quad \sum_{n=1}^{\infty}\left[\psi(1+nz) - \log(nz) - \frac{1}{2nz}\right] + \frac{1}{2z}[\gamma - \log(2\pi z)]$$

$$= \frac{1}{z}\sum_{n=1}^{\infty}\left[\psi\left(1+\frac{n}{z}\right) - \log\left(\frac{n}{z}\right) - \frac{z}{2n}\right] + \frac{1}{2}\left[\gamma - \log\left(\frac{2\pi}{z}\right)\right]$$

and with $z = 2$ we have

$$\sum_{n=1}^{\infty}\left[\psi(1+2n) - \log(2n) - \frac{1}{4n}\right] + \frac{1}{4}[\gamma - \log(4\pi)]$$

$$= \frac{1}{2}\sum_{n=1}^{\infty}\left[\psi\left(1+\frac{n}{2}\right) - \log\left(\frac{n}{2}\right) - \frac{1}{n}\right] + \frac{1}{2}[\gamma - \log \pi]$$

□

We recall (3.4.29)

$$\log G(1+x) = \frac{1}{2}x\log(2\pi) - \frac{1}{2}x(1+x) - \frac{1}{2}\gamma x^2 + \sum_{n=1}^{\infty}\left[\frac{1}{2n}x^2 - x + n\log\left(1+\frac{x}{n}\right)\right]$$

Batir [8] noted that

$$\frac{1}{2n}x^2 - x + n\log\left(1+\frac{x}{n}\right) = \int_0^x \left(\frac{1}{k} - \frac{1}{t+k}\right)t\, dt$$

and we have the summation

$$\sum_{n=1}^{\infty}\left[\frac{1}{2n}x^2 - x + n\log\left(1+\frac{x}{n}\right)\right] = \sum_{n=1}^{\infty}\int_0^x \left(\frac{1}{k} - \frac{1}{t+k}\right) t\, dt$$

$$= \int_0^x \sum_{n=1}^{\infty}\left(\frac{1}{k} - \frac{1}{t+k}\right) t\, dt$$

$$= \int_0^x [\psi(1+t)+\gamma]\, t\, dt$$

where we have used (3.4.14).

Therefore having regard to (3.4.29) we see that

$$\log G(1+x) = \frac{1}{2}x\log(2\pi) - \frac{1}{2}x(1+x) - \frac{1}{2}\gamma x^2 + \int_0^x [\psi(1+t)+\gamma]\, t\, dt$$

and integration by parts results in an easy derivation of Alexeiewsky's theorem [49, p.32]

(3.4.36) $\quad \log G(1+x) = \frac{1}{2}x\log(2\pi) - \frac{1}{2}x(1+x) + x\log\Gamma(1+x) - \int_0^x \log\Gamma(1+t)\, dt$

$\square$

Letting $f(x) = \dfrac{\sin ax}{x}$ in (2.15) gives us

$$\frac{1}{2}a + \sum_{n=1}^{\infty}\frac{\sin na}{n} = \int_0^{\infty}\frac{\sin ax}{x}\, dx + 2\sum_{n=1}^{\infty}\int_0^{\infty}\frac{\sin ax \cos 2\pi nx}{x}\, dx$$

where we have used $f(0) = a$.

It is well known that

$$\int_0^{\infty}\frac{\sin ax}{x}\, dx = \int_0^{\infty}\frac{\sin v}{v}\, dv = \frac{\pi}{2}$$

and we have

$$\int\limits_0^N \frac{\sin ax \cos 2\pi nx}{x} dx = \frac{1}{2}\int\limits_0^N \frac{\sin(2\pi n + a)x - \sin(2\pi n - a)x}{x} dx$$

$$= \frac{1}{2}[Si(2\pi n + a)x - Si(2\pi n - a)x]\Big|_0^N$$

$$= \frac{1}{2}[Si(2\pi n + a)N - Si(2\pi n - a)N]$$

$$= \frac{1}{2}\left[\int\limits_0^{(2\pi n+a)N} \frac{\sin x}{x} dx - \int\limits_0^{(2\pi n-a)N} \frac{\sin ax}{x} dx\right]$$

$$= \frac{1}{2}\int\limits_{(2\pi n-a)N}^{(2\pi n+a)N} \frac{\sin x}{x} dx$$

which gives us the limit

$$\int\limits_0^\infty \frac{\sin ax \cos 2\pi nx}{x} dx = 0$$

We therefore obtain the well known Fourier series

$$\sum_{n=1}^\infty \frac{\sin na}{n} = \frac{1}{2}(\pi - a)$$

□

Letting $f(x) = \dfrac{\cos ax - \cos \pi x}{x}$ in (2.15) gives us

$$\sum_{n=1}^\infty \frac{\cos na - (-1)^n}{n} = \int\limits_0^\infty \frac{\cos ax - \cos \pi x}{x} dx + 2\sum_{n=1}^\infty \int\limits_0^\infty \frac{[\cos ax - \cos \pi x]\cos 2\pi nx}{x} dx$$

where we have used $f(0) = 0$.

$$\int\limits_0^\mu \frac{\cos ax - \cos bx}{x} dx = \int\limits_0^\mu \frac{\cos ax - 1 + 1 - \cos bx}{x} dx$$

$$= \int\limits_0^\mu \frac{\cos ax - 1}{x} dx - \int\limits_0^\mu \frac{\cos bx - 1}{x} dx$$

$$= \int_0^{a\mu} \frac{\cos x - 1}{x} dx - \int_0^{b\mu} \frac{\cos x - 1}{x} dx$$

$$= Ci(a\mu) - Ci(b\mu) + \log \frac{b}{a}$$

where we have used (3.1.3). We therefore obtain the limit

$$\int_0^\infty \frac{\cos ax - \cos bx}{x} dx = \log \frac{b}{a}$$

as reported in [5, p.301] and [48, p.282]. In particular we have

$$\int_0^\infty \frac{\cos ax - \cos \pi x}{x} dx = \log \frac{\pi}{a}$$

We have

$$\int_0^\infty \frac{[\cos ax - \cos \pi x]\cos 2\pi nx}{x} dx = \frac{1}{2}\int_0^\infty \frac{\cos(2\pi n + a)x + \cos(2\pi n - a)x}{x} dx$$

$$- \frac{1}{2}\int_0^\infty \frac{\cos(2n+1)\pi x + \cos(2n-1)\pi x}{x} dx$$

$$= \frac{1}{2}\int_0^\infty \frac{\cos(2\pi n + a)x - \cos(2n+1)\pi x}{x} dx + \frac{1}{2}\int_0^\infty \frac{\cos(2\pi n - a)x - \cos(2n-1)\pi x}{x} dx$$

$$= \frac{1}{2}\log \frac{(2n+1)\pi}{(2\pi n + a)} + \frac{1}{2}\log \frac{(2n-1)\pi}{(2\pi n - a)}$$

We see that

$$\sum_{n=1}^\infty \frac{\cos na - (-1)^n}{n} = \sum_{n=1}^\infty \frac{\cos na}{n} + \log 2$$

and using the familiar trigonometric series shown in Carslaw's book [17, p.241]

$$\log[2\sin(a/2)] = -\sum_{n=1}^\infty \frac{\cos na}{n} \quad (0 < a < 2\pi)$$

we obtain

$$\log \sin(a\pi/2) = \log a + \sum_{n=1}^{\infty} \log \frac{4n^2 - a^2}{4n^2 - 1}$$

$$= \log a + \sum_{n=1}^{\infty} \log \frac{4n^2 - a^2}{4n^2} - \sum_{n=1}^{\infty} \log \frac{4n^2 - 1}{4n^2}$$

Using the Wallis product

$$\frac{\pi}{2} = \lim_{N \to \infty} \prod_{n=1}^{N} \frac{2n}{2n-1} \cdot \frac{2n}{2n+1}$$

we see that

$$\log \frac{\pi}{2} = -\sum_{n=1}^{\infty} \log \frac{4n^2 - 1}{4n^2}$$

and we obtain

$$\log \sin(\pi a/2) = \log(\pi a/2) + \sum_{n=1}^{\infty} \log \frac{4n^2 - a^2}{4n^2}$$

or equivalently

$$\log \sin(\pi a) = \log(\pi a) + \sum_{n=1}^{\infty} \log \frac{n^2 - a^2}{n^2}$$

which corresponds with the Euler product formula for the sine function

$$\sin \pi a = \pi a \prod_{n=1}^{\infty} \left(1 - \frac{a^2}{n^2}\right)$$

$\square$

Integration of (3.4.7) with respect to $t$ results in

$$(3.4.7) \int_0^{\infty} \left[ \log \Gamma(x+u) - \log \Gamma(u) - (x+u) \log(x+u) + u \log u + x + \frac{1}{2} \log(x+u) - \frac{1}{2} \log u \right] du$$

$$= \frac{x}{4}[-x + 2(x-1) \log x] + \log G(1+x) - x \log \Gamma(x)$$

where we have used (3.4.36).

Integrating (3.4.20) gives us

$$\sum_{n=1}^{\infty}\left[\log\Gamma(x+n)-\log\Gamma(n)-x\psi(1+n)-\frac{1}{1+n}\left(\frac{1}{2}x^2-x\right)\right]=-\frac{1}{2}x^2\gamma+2\left(\frac{1}{2}x^2-x\right)-\int_0^x u\psi(u)du$$

and we have

$$\int_0^x u\psi(u)du = x\log\Gamma(x)-\int_0^x \log\Gamma(u)\,du$$

This gives us

$$\sum_{n=1}^{\infty}\left[\log\Gamma(x+n)-\log\Gamma(n)-x\psi(1+n)-\frac{1}{1+n}\left(\frac{1}{2}x^2-x\right)\right]=-\frac{1}{2}\gamma x^2+\frac{1}{2}(x^2-3x)$$

$$-\log G(1+x)+\frac{1}{2}x\log(2\pi)$$

It is easily seen that this concurs with (3.4.36).

### 3.5 Stieltjes constants

The Stieltjes constants $\gamma_p(x)$ are the coefficients of the Laurent expansion of the Hurwitz zeta function $\varsigma(s,x)$ about $s=1$

(3.5.1) $$\varsigma(s,x)=\sum_{n=0}^{\infty}\frac{1}{(n+x)^s}=\frac{1}{s-1}+\sum_{p=0}^{\infty}\frac{(-1)^p}{p!}\gamma_p(x)(s-1)^p$$

and we have [54]

(3.5.2) $$\gamma_0(x)=-\psi(x)$$

With $x=1$ equation (3.5.1) reduces to the Riemann zeta function

$$\varsigma(s)=\frac{1}{s-1}+\sum_{p=0}^{\infty}\frac{(-1)^p}{p!}\gamma_p(s-1)^p$$

where $\gamma_p(1)=\gamma_p$.

We write (2.15) in the form

(3.5.2.1) $$\frac{1}{2}f(0) + \lim_{N\to\infty}\left[\sum_{n=1}^{N} f(n) - \int_0^N f(x)\,dx\right] = 2\sum_{n=1}^{\infty}\int_0^{\infty} f(x)\cos 2\pi nx\,dx$$

and, with the function $f(x) = \dfrac{\log^p(1+x)}{1+x}$, this results in

$$\lim_{N\to\infty}\left[\sum_{n=1}^{N}\frac{\log^p(1+n)}{1+n} - \frac{\log^{p+1}(1+N)}{1+p}\right] = 2\sum_{n=1}^{\infty}\int_0^{\infty}\frac{\log^p(1+x)}{1+x}\cos 2\pi nx\,dx$$

since $f(0) = 0$ for $p \geq 1$.

Simple algebra gives us

$$\sum_{n=1}^{N}\frac{\log^p(1+n)}{1+n} - \frac{\log^{p+1}(1+N)}{1+p} = \sum_{n=1}^{N}\frac{\log^p n}{n} + \frac{\log^p(1+N)}{1+N} - \frac{\log^{p+1}(1+N)}{1+p}$$

$$= \sum_{n=1}^{N}\frac{\log^p n}{n} - \frac{\log^{p+1} N}{1+p} + \left[\frac{\log^{p+1} N}{1+p} - \frac{\log^{p+1}(1+N)}{1+p}\right] + \frac{\log^p(1+N)}{1+N}$$

and it is easily shown that the expression in parentheses vanishes as $N\to\infty$. Successive applications of L'Hôpital's rule shows that $\displaystyle\lim_{N\to\infty}\frac{\log^p(1+N)}{1+N} = 0$ and hence we have

(3.5.3) $$\lim_{N\to\infty}\left[\sum_{n=1}^{N}\frac{\log^p(1+n)}{1+n} - \frac{\log^{p+1}(1+N)}{1+p}\right] = \lim_{N\to\infty}\left[\sum_{n=1}^{N}\frac{\log^p n}{n} - \frac{\log^{p+1} N}{1+p}\right]$$

It is well known that Stieltjes [35] proved in 1885 that the Stieltjes constants $\gamma_p$ may be represented by

(3.5.4) $$\lim_{N\to\infty}\left[\sum_{n=1}^{N}\frac{\log^p n}{n} - \frac{\log^{p+1} N}{1+p}\right] = \gamma_p$$

We see that

$$\int_0^{\infty}\cos 2\pi nx\,\frac{\log^p(1+x)}{1+x}\,dx = \int_0^{\infty}\cos 2\pi n(1+x)\,\frac{\log^p(1+x)}{1+x}\,dx$$

and an elementary substitution gives us

$$\int_0^\infty \cos 2\pi n(1+x) \frac{\log^p(1+x)}{1+x} dx = \int_1^\infty \cos 2\pi nu \frac{\log^p u}{u} du$$

Hence we obtain for $p \geq 1$

(3.5.5) $$\gamma_p = 2\sum_{n=1}^\infty \int_1^\infty \cos 2\pi nx \frac{\log^p x}{x} dx$$

which was originally derived by Briggs [14] in 1955.

It should be noted that this only applies for $p \geq 1$. From Knopp's book [38, p.521] the Euler-Maclaurin summation formula gives us

(3.5.6) $$\gamma_0 = \gamma = \frac{1}{2} + \sum_{n=1}^\infty \int_1^\infty \frac{\sin 2\pi nx}{n\pi x^2} dx$$

and integration by parts gives us

$$\int_1^\infty \frac{\cos 2\pi nx}{x} dx = \frac{\sin 2\pi nx}{2\pi nx}\Big|_1^\infty + \int_1^\infty \frac{\sin 2\pi nx}{2\pi nx^2} dx = \int_1^\infty \frac{\sin 2\pi nx}{2\pi nx^2} dx$$

Hence we have an additional factor of 1/2 in the case where $p = 0$

(3.5.7) $$\gamma_0 = \gamma = \frac{1}{2} + 2\sum_{n=1}^\infty \int_1^\infty \frac{\cos 2\pi nx}{x} dx$$

This factor arises because $f(0) = 1$ in (3.5.2.1) when $p = 0$.

We note from (3.5.5) that

$$\gamma_p = 2(-1)^p \frac{d^p}{ds^p} \sum_{n=1}^\infty \int_1^\infty \frac{\cos 2\pi nx}{x^s} dx \Big|_{s=1}$$

□

The above analysis may be generalised by considering the function $f(x) = \frac{\log^p(a+x)}{a+x}$. Proceeding as before results in

$$\frac{1}{2}\frac{\log^p a}{a} + \lim_{N \to \infty}\left[\sum_{n=1}^N \frac{\log^p(a+n)}{a+n} - \frac{\log^{p+1}(a+N) - \log^{p+1} a}{p+1}\right] = 2\sum_{n=1}^\infty \int_0^\infty \frac{\log^p(a+x)}{a+x} \cos 2\pi nx\, dx$$

which may be written as

$$-\frac{1}{2}\frac{\log^p a}{a}+\frac{\log^{p+1} a}{p+1}+\lim_{N\to\infty}\left[\sum_{n=0}^{N}\frac{\log^p(a+n)}{a+n}-\frac{\log^{p+1}(a+N)}{p+1}\right]=2\sum_{n=1}^{\infty}\int_0^{\infty}\frac{\log^p(a+x)}{a+x}\cos 2\pi nx\,dx$$

We will see below in (3.7.17) that for $p \geq 1$

$$\gamma_p(a)=\frac{1}{2}\frac{\log^p a}{a}-\frac{\log^{p+1} a}{p+1}+2\sum_{n=1}^{\infty}\int_0^{\infty}\frac{\log^p(a+x)}{a+x}\cos 2\pi nx\,dx$$

and hence we have

(3.5.8) $$\gamma_p(a)=\lim_{N\to\infty}\left[\sum_{n=0}^{N}\frac{\log^p(a+n)}{a+n}-\frac{\log^{p+1}(a+N)}{p+1}\right]$$

as previously shown by Berndt [9] in 1972. With $a = 1$ in (3.5.8) and using (3.5.3) we see that Stieltjes's formula (3.5.4) immediately follows.

□

We see that

$$\gamma_p(1+a)=\frac{1}{2}\frac{\log^p(1+a)}{1+a}-\frac{\log^{p+1}(1+a)}{p+1}+2\sum_{n=1}^{\infty}\int_0^{\infty}\frac{\log^p(a+1+x)}{a+1+x}\cos 2\pi nx\,dx$$

$$=\frac{1}{2}\frac{\log^p(1+a)}{1+a}-\frac{\log^{p+1}(1+a)}{p+1}+2\sum_{n=1}^{\infty}\int_0^{\infty}\frac{\log^p(a+1+x)}{a+1+x}\cos 2\pi n(1+x)\,dx$$

$$=\frac{1}{2}\frac{\log^p(1+a)}{1+a}-\frac{\log^{p+1}(1+a)}{p+1}+2\sum_{n=1}^{\infty}\int_1^{\infty}\frac{\log^p(a+x)}{a+x}\cos 2\pi nx\,dx$$

Therefore we have

$$\gamma_p(1+a)-\gamma_p(a)$$

$$=\frac{1}{2}\frac{\log^p(1+a)}{1+a}-\frac{1}{2}\frac{\log^p a}{a}-\frac{\log^{p+1}(1+a)}{p+1}+\frac{\log^{p+1} a}{p+1}-2\sum_{n=1}^{\infty}\int_0^{1}\frac{\log^p(a+x)}{a+x}\cos 2\pi nx\,dx$$

We see from (2.13) that

$$2\sum_{n=1}^{\infty}\int_0^1 \frac{\log^p(a+x)}{a+x}\cos 2\pi nx\, dx = \frac{1}{2}\left[\frac{\log^p a}{a}+\frac{\log^p(a+1)}{a+1}\right]-\int_0^1 \frac{\log^p(a+x)}{a+x}dx$$

$$=\frac{1}{2}\left[\frac{\log^p a}{a}+\frac{\log^p(a+1)}{a+1}\right]-\left[\frac{\log^{p+1}(a+1)}{p+1}-\frac{\log^{p+1} a}{p+1}\right]$$

and so we obtain

(3.5.9) $$\gamma_p(1+a)-\gamma_p(a)=-\frac{\log^p a}{a}$$

This identity may be obtained more directly as follows:

Since $\varsigma(s,1+a)-\varsigma(s,a)=-\frac{1}{a^s}$ we have

$$\varsigma^{(p)}(s,1+a)-\varsigma^{(p)}(s,a)=(-1)^{p+1}\frac{\log^p a}{a^s}$$

and therefore using

$$\gamma_p(x)-\gamma_p(y)=\lim_{s\to 1}(-1)^p\frac{\partial^p}{\partial s^p}[\varsigma(s,x)-\varsigma(s,y)]$$

we obtain

$$\gamma_p(1+a)-\gamma_p(a)=-\frac{\log^p a}{a}$$

**3.6 Connections with other summation formulae**

The following integral formula was originally obtained by Coffey in 2007 for the Stieltjes constants

(3.6.1) $$\gamma_p(a)=\frac{1}{2}\frac{\log^p a}{a}-\frac{\log^{p+1} a}{p+1}+i\int_0^\infty \frac{(a-ix)\log^p(a+ix)-(a+ix)\log^p(a-ix)}{(a^2+x^2)(e^{2\pi x}-1)}dx$$

and a different derivation appears in [24].

It should be noted that

$$\int_0^\infty \frac{(a-ix)\log^p(a+ix)-(a+ix)\log^p(a-ix)}{(a^2+x^2)(e^{2\pi x}-1)}dx$$

$$= \int_0^\infty \left[ \frac{\log^p(a+ix)}{a+ix} - \frac{\log^p(a-ix)}{a-ix} \right] \frac{dx}{e^{2\pi x}-1}$$

The structure of the integral in (3.6.1) therefore indicates the close connection with the Abel-Plana summation formula [49, p.90]

(3.6.2) $$\frac{1}{2}f(0) + \sum_{n=1}^\infty f(n) = \int_0^\infty f(x)\,dx + i\int_0^\infty \frac{f(ix)-f(-ix)}{e^{2\pi x}-1}\,dx$$

which applies to functions which are analytic in the right-hand plane and satisfy the convergence condition $\lim_{y\to\infty} e^{-2\pi y}|f(x+iy)| = 0$ uniformly on any finite interval of $x$. Derivations of the Abel-Plana summation formula are given in [6], [34, p.339], [52, p.145] and [51, p.118].

Adamchik [4] noted that the Hermite integral for the Hurwitz zeta function may be derived from the Abel-Plana summation formula.

□

With $f(x) \to f(2x)$ in (3.6.2) we have

(3.6.3) $$\frac{1}{2}f(0) + \sum_{n=1}^\infty f(2n) = \int_0^\infty f(2x)\,dx + i\int_0^\infty \frac{f(2ix)-f(-2ix)}{e^{2\pi x}-1}\,dx$$

$$= \frac{1}{2}\int_0^\infty f(u)\,du + \frac{1}{2}i\int_0^\infty \frac{f(iu)-f(-iu)}{e^{\pi u}-1}\,du$$

Subtraction results in

$$\frac{1}{2}f(0) + 2\sum_{n=1}^\infty f(2n) - \sum_{n=1}^\infty f(n) = i\int_0^\infty \frac{f(ix)-f(-ix)}{e^{\pi x}-1}\,dx - i\int_0^\infty \frac{f(ix)-f(-ix)}{e^{2\pi x}-1}\,dx$$

and, since for suitably convergent series $2\sum_{n=1}^\infty a_{2n} - \sum_{n=1}^\infty a_n = \sum_{n=1}^\infty (-1)^n a_n$, we have

(3.6.4) $$\frac{1}{2}f(0) + \sum_{n=1}^\infty (-1)^n f(n) = i\int_0^\infty \frac{f(ix)-f(-ix)}{2\sinh \pi x}\,dx$$

which is the alternating form of the Abel-Plana summation formula originally derived by Abel [16]. This may be compared with (2.16).

For completeness, we note that Saharian [46] reports the summation

(3.6.5) $$\sum_{n=0}^{\infty} f\left(n+\frac{1}{2}\right) = \int_0^{\infty} f(x)\,dx - i\int_0^{\infty} \frac{f(ix)-f(-ix)}{e^{2\pi x}+1}\,dx$$

□

Letting $f(x) \to f(x)\cos \pi x$ in (3.6.2) results in

$$\frac{1}{2}f(0) + \sum_{n=1}^{\infty}(-1)^n f(n) = \int_0^{\infty} f(x)\cos \pi x\,dx + i\int_0^{\infty} \frac{f(ix)-f(-ix)}{e^{2\pi x}-1}\cos i\pi x\,dx$$

and this may be written as

(3.6.6) $$\frac{1}{2}f(0) + \sum_{n=1}^{\infty}(-1)^n f(n) = \int_0^{\infty} f(x)\cos \pi x\,dx + i\int_0^{\infty} \frac{f(ix)-f(-ix)}{e^{2\pi x}-1}\cosh \pi x\,dx$$

and equating this with (3.6.4) gives us

$$i\int_0^{\infty} \frac{f(ix)-f(-ix)}{2\sinh \pi x}\,dx = \int_0^{\infty} f(x)\cos \pi x\,dx + i\int_0^{\infty} \frac{f(ix)-f(-ix)}{e^{2\pi x}-1}\cosh \pi x\,dx$$

which simplifies to

$$\int_0^{\infty} f(x)\cos \pi x\,dx = -\frac{1}{2}i\int_0^{\infty} \frac{[f(ix)-f(-ix)]e^{-\pi x}}{e^{2\pi x}-1}\,dx$$

□

Comparing (3.6.2) and (2.15) we see that

(3.6.7) $$i\int_0^{\infty} \frac{f(ix)-f(-ix)}{e^{2\pi x}-1}\,dx = 2\sum_{n=1}^{\infty}\int_0^{\infty} f(x)\cos 2\pi nx\,dx$$

for functions which satisfy the conditions relevant to both the Abel-Plana summation formula and the Poisson summation formula. An application of (3.6.7) is shown below.

With $f(x) = e^{-ax}$ we have

$$i\int_0^{\infty} \frac{f(ix)-f(-ix)}{e^{2\pi x}-1}\,dx = 2\int_0^{\infty} \frac{\sin(ax)}{e^{2\pi x}-1}\,dx$$

Therefore (3.6.7) gives us

$$\int_0^\infty \frac{\sin(ax)}{e^{2\pi x}-1}dx = \sum_{n=1}^\infty \int_0^\infty e^{-ax}\cos 2\pi nx\, dx$$

We easily obtain the indefinite integral

$$\int e^{-ax}\cos 2\pi nx\, dx = \frac{e^{-ax}[2\pi n\sin 2\pi n - a\cos 2\pi nx]}{a^2+4\pi^2 n^2}$$

so that

(3.6.8) $$\int_0^\infty e^{-ax}\cos 2\pi nx\, dx = \frac{a}{a^2+4\pi^2 n^2}$$

We then have

$$\int_0^\infty \frac{\sin(ax)}{e^{2\pi x}-1}dx = \sum_{n=1}^\infty \frac{a}{a^2+4\pi^2 n^2}$$

We may then use the well known identity ([15, p.296], [38, p.378])

(3.6.9) $$2\sum_{n=1}^\infty \frac{a}{a^2+4\pi^2 n^2} = \frac{1}{e^a-1} - \frac{1}{a} + \frac{1}{2}$$

to determine Legendre's relation [52, p.119]

(3.6.10) $$2\int_0^\infty \frac{\sin(ax)}{e^{2\pi x}-1}dx = \frac{1}{e^a-1} - \frac{1}{a} + \frac{1}{2} = \frac{1}{2}\coth\frac{a}{2} - \frac{1}{a}$$

A rigorous derivation of this result is shown in Bromwich's book [15, p.501]. It is interesting to note that Apostol [5, p.334] employed Poisson's formula to derive the corresponding partial fraction decomposition formula for $\coth x$. A slightly simplified version of this proof is set out below.

With $f(x) = e^{-ax}$ in (2.15) we have

$$\frac{1}{2} + \sum_{n=1}^\infty e^{-nx} = \int_0^\infty e^{-ax}dx + 2\sum_{n=1}^\infty \int_0^\infty e^{-ax}\cos 2\pi nx\, dx$$

so that upon using (3.6.8) we immediately obtain (3.6.9).

Various applications of the Abel-Plana summation formula are contained in Ramanujan's Notebooks, Part V [11].

As noted by Ivić [36, p.6] the identity (3.6.9) may be readily obtained from Euler's infinite product representation of the $\sin x$ function

$$\sin x = x \prod_{n=1}^{\infty} \left(1 - \frac{x^2}{n^2 \pi^2}\right)$$

which was derived earlier in this paper using the Poisson summation formula.

The substitution $x = -\frac{1}{2} iu$ gives us

$$\sin\left(\frac{1}{2} iu\right) = \frac{1}{2} iu \prod_{n=1}^{\infty} \left(1 + \frac{u^2}{(2\pi n)^2}\right)$$

so that

$$\sinh\left(\frac{1}{2} u\right) = \frac{1}{2} u \prod_{n=1}^{\infty} \left(\frac{(2\pi n)^2 + u^2}{(2\pi n)^2}\right)$$

Logarithmic differentiation then results in

$$\frac{1}{2} \coth\left(\frac{1}{2} u\right) = \frac{1}{u} + 2 \sum_{n=1}^{\infty} \frac{u}{u^2 + 4\pi^2 n^2}$$

which is equivalent to (3.6.9).

□

Letting $f(x) = \frac{1}{u+x}$ in (3.6.7) gives us

$$\int_0^{\infty} \frac{x}{(u^2 + x^2)(e^{2\pi x} - 1)} dx = \sum_{n=1}^{\infty} \int_0^{\infty} \frac{\cos 2\pi nx}{u + x} dx$$

$$= \sum_{n=1}^{\infty} [\cos(2n\pi u) Ci(2n\pi u) + \sin(2n\pi u) si(2n\pi u)]$$

and using (3.1.7)

$$\psi(u) = \log u - \frac{1}{2u} + 2 \sum_{n=1}^{\infty} [\cos(2n\pi u) Ci(2n\pi u) + \sin(2n\pi u) si(2n\pi u)]$$

we obtain Binet's formula (3.4.5)

$$\log u - \psi(u) = \frac{1}{2u} + 2\int_0^\infty \frac{x}{(u^2 + x^2)(e^{2\pi x} - 1)} dx$$

□

We also note the recent paper by Butzer et al. [16] "The Summation Formulae of Euler–Maclaurin, Abel–Plana, Poisson, and their interconnections with the Approximate Sampling Formula of Signal Analysis" where it is shown that these four fundamental formulae are all equivalent, in the sense that each is a corollary of any of the others.

As mentioned in [16], the Abel-Plana summation formula may be expressed as

$$(3.6.11) \qquad \frac{1}{2} f(0) + \lim_{N \to \infty} \left[ \sum_{n=1}^N f(n) - \int_0^{N+\frac{1}{2}} f(x) dx \right] = i \int_0^\infty \frac{f(ix) - f(-ix)}{e^{2\pi x} - 1} dx$$

and with $f(x) = x \log(1+x)$ we have

$$i \int_0^\infty \frac{f(ix) - f(-ix)}{e^{2\pi x} - 1} dx = -\int_0^\infty \frac{x \log(1+x^2)}{e^{2\pi x} - 1} dx$$

This gives us

$$(3.6.12) \qquad \lim_{N \to \infty} \left[ \int_0^{N+\frac{1}{2}} x \log(1+x) dx - \sum_{n=1}^N n \log(1+n) \right] = \int_0^\infty \frac{x \log(1+x^2)}{e^{2\pi x} - 1} dx$$

We note that

$$\sum_{n=1}^N n \log(1+n) = \sum_{m=2}^{N+1} (m-1) \log m$$

$$= \sum_{m=1}^{N+1} (m-1) \log m$$

$$= \sum_{n=1}^N (n-1) \log n + N \log(N+1)$$

$$= \sum_{n=1}^N (n-1) \log n + N \log N + N \log\left(1 + \frac{1}{N}\right)$$

and we have

$$\int_0^{N+\frac{1}{2}} x\log(1+x)\,dx = \frac{1}{2}\left(N^2 + N - \frac{3}{4}\right)\log\left(N+\frac{3}{2}\right) - \frac{1}{4}\left(N^2 + N + \frac{1}{4}\right) + \frac{1}{2}\left(N+\frac{1}{2}\right)$$

$$= \frac{1}{2}\left(N^2 + N - \frac{3}{4}\right)\log N + \frac{1}{2}\left(N^2 + N - \frac{3}{4}\right)\log\left(1+\frac{3}{2N}\right) - \frac{1}{4}\left(N^2 + N + \frac{1}{4}\right) + \frac{1}{2}\left(N+\frac{1}{2}\right)$$

Designating $I_N = \int_0^{N+\frac{1}{2}} x\log(1+x)\,dx - \sum_{n=1}^{N} n\log(1+n)$ we then have

$$I_N = \int_0^{N+\frac{1}{2}} x\log(1+x)\,dx - \sum_{n=1}^{N} n\log(1+n) = -\sum_{n=1}^{N}(n-1)\log n - N\log\left(1+\frac{1}{N}\right)$$

$$+ \frac{1}{2}\left(N^2 - N - \frac{3}{4}\right)\log N + \frac{1}{2}\left(N^2 + N - \frac{3}{4}\right)\log\left(1+\frac{3}{2N}\right) - \frac{1}{4}\left(N^2 + N + \frac{1}{4}\right) + \frac{1}{2}\left(N+\frac{1}{2}\right)$$

Since $\log\left(1+\frac{u}{N}\right) = \frac{u}{N} - \frac{u^2}{2N^2} + O(N^{-3})$ we have

$$\frac{1}{2}\left(N^2 + N - \frac{3}{4}\right)\log\left(1+\frac{3}{2N}\right) = \frac{1}{2}\left(N^2 + N - \frac{3}{4}\right)\left[\frac{3}{2N} - \frac{9}{8N^2} + O(N^{-3})\right]$$

$$= \frac{3}{4}N - \frac{9}{16} + \frac{3}{4} + O(N^{-1})$$

We then have

(3.6.13) $$\lim_{N\to\infty}\left[\int_0^{N+\frac{1}{2}} x\log(1+x)\,dx - \sum_{n=1}^{N} n\log(1+n)\right]$$

$$= \lim_{N\to\infty}\left[-\sum_{n=1}^{N}(n-1)\log n + \frac{1}{2}\left(N^2 - N - \frac{3}{4}\right)\log N + N - \frac{1}{4}N^2 - \frac{5}{8}\right]$$

Using the Euler-Maclaurin summation formula, Hardy [34, p.333] showed that the Riemann zeta function could be expressed as follows:

$$\varsigma(s) = \lim_{n\to\infty}\left[\sum_{k=1}^{n}\frac{1}{k^s} - \frac{n^{1-s}}{1-s} - \frac{1}{2}n^{-s}\right] \qquad \text{Re}(s) > -1$$

$$\varsigma(s) = \lim_{n\to\infty}\left[\sum_{k=1}^{n}\frac{1}{k^s} - \frac{n^{1-s}}{1-s} - \frac{1}{2}n^{-s} + \frac{1}{12}sn^{-s-1}\right] \qquad \text{Re}(s) > -3$$

Differentiating the first identity results in for $\text{Re}(s) > -1$

$$\varsigma'(s) = \lim_{n\to\infty}\left[-\sum_{k=1}^{n}\frac{\log k}{k^s} + \frac{n^{1-s}(1-s)\log n - n^{1-s}}{(1-s)^2} + \frac{1}{2}n^{-s}\log n\right]$$

and with $s = 0$ we obtain

$$\varsigma'(0) = \lim_{n\to\infty}\left[-\sum_{k=1}^{n}\log k + \left(n+\frac{1}{2}\right)\log n - n\right]$$

Hence, using the Stirling approximation we see that $\varsigma'(0) = -\frac{1}{2}\log(2\pi)$.

For $\text{Re}(s) > -3$ we have

$$\varsigma'(s) = \lim_{n\to\infty}\left[-\sum_{k=1}^{n}\frac{\log k}{k^s} + \frac{n^{1-s}(1-s)\log n - n^{1-s}}{(1-s)^2} + \frac{1}{2}n^{-s}\log n - \frac{1}{12}sn^{-s-1}\log n + \frac{1}{12}n^{-s-1}\right]$$

and with $s = -1$ we get (cf [49, p.25] re the definition of the Glaisher-Kinkelin constant)

$$\varsigma'(-1) = \lim_{n\to\infty}\left[-\sum_{k=1}^{n}k\log k + \left(\frac{n^2}{2}+\frac{n}{2}+\frac{1}{12}\right)\log n - \frac{n^2}{4} + \frac{1}{12}\right]$$

We then have

$$(3.6.14) \quad \varsigma'(-1) - \frac{3}{4} + \frac{1}{2}\log(2\pi) = \lim_{n\to\infty}\left[-\sum_{k=1}^{n}(k-1)\log k + \left(\frac{n^2}{2}-\frac{n}{2}-\frac{5}{12}\right)\log n + n - \frac{n^2}{4} + \frac{2}{3}\right]$$

Adamchik [3] has shown that

$$(3.6.15) \qquad \int_0^{\infty}\frac{x\log(1+x^2)}{e^{2\pi x}-1}dx = \varsigma'(-1) - \frac{3}{4} + \frac{1}{2}\log(2\pi)$$

It is clear that further work is required to resolve the two discrepancies between (3.6.13) and (3.6.14).

**3.7 Hurwitz zeta function**

We now consider the function $f(x) = \dfrac{1}{(a+x)^s}$, where $\operatorname{Re}(s) > 1$, and using (2.15) we obtain

$$\frac{1}{2}\frac{1}{a^s} + \sum_{n=1}^{\infty} \frac{1}{(a+n)^s} = \int_0^{\infty} \frac{dx}{(a+x)^s} + 2\sum_{n=1}^{\infty} \int_0^{\infty} \frac{\cos 2\pi nx}{(a+x)^s}\,dx$$

which may be written as

$$\sum_{n=0}^{\infty} \frac{1}{(a+n)^s} = \frac{1}{2}\frac{1}{a^s} + \frac{1}{(s-1)a^{s-1}} + 2\sum_{n=1}^{\infty} \int_0^{\infty} \frac{\cos 2\pi nx}{(a+x)^s}\,dx$$

The Hurwitz zeta function is defined for $\operatorname{Re}(s) > 1$ as

$$\varsigma(s,a) = \sum_{n=0}^{\infty} \frac{1}{(a+n)^s}$$

and we therefore have

(3.7.1) $$\varsigma(s,a) = \frac{1}{2}\frac{1}{a^s} + \frac{1}{(s-1)a^{s-1}} + 2\sum_{n=1}^{\infty} \int_0^{\infty} \frac{\cos 2\pi nx}{(a+x)^s}\,dx$$

which corrects a misprint in Mordell's paper [41].

Using (2.13) we obtain

$$\frac{1}{2}\left[\frac{1}{a^s} + \frac{1}{(a+1)^s}\right] = \int_0^1 \frac{dx}{(a+x)^s} + 2\sum_{n=1}^{\infty} \int_0^1 \frac{\cos 2\pi nx}{(a+x)^s}\,dx$$

and hence we have

(3.7.2) $$\frac{1}{2}\left[\frac{1}{a^s} + \frac{1}{(a+1)^s}\right] = -\frac{1}{s-1}\left[\frac{1}{(a+1)^{s-1}} - \frac{1}{a^{s-1}}\right] + 2\sum_{n=1}^{\infty} \int_0^1 \frac{\cos 2\pi nx}{(a+x)^s}\,dx$$

Subtracting (3.7.2) from (3.7.1) gives us

(3.7.3) $$\varsigma(s,a) = \frac{1}{a^s} + \frac{1}{2}\frac{1}{(a+1)^s} + \frac{1}{(s-1)(a+1)^{s-1}} + 2\sum_{n=1}^{\infty} \int_1^{\infty} \frac{\cos 2\pi nx}{(a+x)^s}\,dx$$

Letting $s = 0$ in (3.7.1) results in the well-known formula

$$\varsigma(0,a) = \frac{1}{2} - a$$

and we are therefore lead to suspect that (3.7.1) is valid for a larger region than just $\text{Re}(s) > 1$; as demonstrated by Mordell [41], by using integration by parts it is easily shown that $\varsigma(s,a)$ exists over all the $s$-plane except for a simple pole at $s = 1$.

In particular, (3.7.1) is valid for $\text{Re}(s) < 0$ and, in view of this, Mordell [41] was able to show how the representation (3.7.1) may be used to derive Hurwitz's Fourier series for $\varsigma(s,a)$; for ease of reference, his proof is outlined below (with the correction of some misprints).

Assuming that $\text{Re}(s) < 0$, in the case where $1 \geq a > 0$ and the interval of integration is $[-a, 0]$, then the left-hand side of (2.12) gives us

$$(a-a)^{-s} + \frac{1}{2}(a+0)^{-s} = \frac{1}{2a^s}$$

so that (2.12) results in

$$\frac{1}{2a^s} = \int_{-a}^{0} \frac{dx}{(a+x)^s} + 2\sum_{n=1}^{\infty} \int_{-a}^{0} \frac{\cos 2\pi nx}{(a+x)^s} dx$$

$$= -\frac{a^{1-s}}{s-1} + 2\sum_{n=1}^{\infty} \int_{-a}^{0} \frac{\cos 2\pi nx}{(a+x)^s} dx$$

Then, upon adding (3.7.1) to this, we obtain for $\text{Re}(s) < 0$

(3.7.4) $$\varsigma(s,a) = 2\sum_{n=1}^{\infty} \int_{-a}^{\infty} \frac{\cos 2\pi nx}{(a+x)^s} dx$$

and using integration by parts this becomes

$$\varsigma(s,a) = \lim_{N \to \infty} \sum_{n=1}^{\infty} \frac{\sin 2\pi nx}{\pi n(a+x)^s} \bigg|_{-a}^{N} + 2s \sum_{n=1}^{\infty} \int_{-a}^{\infty} \frac{\sin 2\pi nx}{2n\pi(a+x)^{s+1}} dx$$

$$= 2s \sum_{n=1}^{\infty} \int_{-a}^{\infty} \frac{\sin 2\pi nx}{2n\pi(a+x)^{s+1}} dx$$

$$= 2s \sum_{n=1}^{\infty} \int_0^{\infty} \frac{\sin 2\pi n(u-a)}{2n\pi u^{s+1}} du$$

and we obtain

(3.7.5) $$\varsigma(s,a) = 2s \sum_{n=1}^{\infty} \left[ \frac{\cos 2\pi na}{2n\pi} \int_0^{\infty} \frac{\sin 2\pi nu}{u^{s+1}} du - \frac{\sin 2\pi na}{2n\pi} \int_0^{\infty} \frac{\cos 2\pi nu}{u^{s+1}} du \right]$$

Using contour integration we have ([50, p.107], [18, p.91]) for $a > 0$ and $1 > p > 0$

(3.7.6) $$\int_0^{\infty} u^{p-1} \cos bu\, du = \frac{\Gamma(p)\cos(\pi p/2)}{b^p}$$

and with $p = -s$ we see that for $0 > s > -1$

$$\int_0^{\infty} \frac{\cos bu}{u^{s+1}} du = \frac{\Gamma(-s)\cos(\pi s/2)}{b^{-s}}$$

$$= -\frac{\Gamma(1-s)\cos(\pi s/2)}{sb^{-s}}$$

Similarly, we have

(3.7.7) $$\int_0^{\infty} u^{p-1} \sin bu\, du = \frac{\Gamma(p)\sin(\pi p/2)}{b^p}$$

so that

$$\int_0^{\infty} \frac{\sin bu}{u^{s+1}} du = \frac{\Gamma(1-s)\sin(\pi s/2)}{sb^{-s}}$$

We therefore obtain Hurwitz's formula

(3.7.8) $$\varsigma(s,a) = 2\Gamma(1-s) \left[ \sin\left(\frac{\pi s}{2}\right) \sum_{n=1}^{\infty} \frac{\cos 2\pi na}{(2\pi n)^{1-s}} + \cos\left(\frac{\pi s}{2}\right) \sum_{n=1}^{\infty} \frac{\sin 2\pi na}{(2\pi n)^{1-s}} \right]$$

where $\mathrm{Re}(s) < 0$ and $0 < a \leq 1$. In 2000, Boudjelkha [13] showed that this formula also applies in the region $\mathrm{Re}(s) < 1$.

With $a = 1$ this reduces to the functional equation for the Riemann zeta function

(3.7.9) $$\varsigma(1-s) = 2(2\pi)^{-s}\Gamma(s)\cos(\pi s/2)\varsigma(s)$$

Employing the same procedure with $f(x) = \dfrac{\cos \pi x}{(a+x)^s}$, it should be possible to derive Boudjelkha's formula [13] for the alternating Hurwitz zeta function $\varsigma_a(s,a)$ which is defined by

$$\varsigma_a(s,a) = \sum_{n=0}^{\infty} \frac{(-1)^n}{(a+n)^s}$$

Boudjelkha's formula is

(3.7.10) $$\varsigma_a(s,a) = 2\Gamma(1-s)\pi^{s-1}\left[\sin\left(\frac{\pi s}{2}\right)\sum_{n=0}^{\infty}\frac{\cos(2n+1)\pi a}{(2n+1)^{1-s}} + \cos\left(\frac{\pi s}{2}\right)\sum_{n=0}^{\infty}\frac{\sin(2n+1)\pi a}{(2n+1)^{1-s}}\right]$$

and holds under the same conditions as (3.7.8) above, namely:

$$(\sigma < 0,\ 0 < a \leq 1;\ 0 < \sigma,\ a < 1)$$

We also note that Oberhettinger [45] used the Poisson summation formula to derive a corresponding formula for the Lerch zeta function $\Phi(z,s,a) = \sum_{n=0}^{\infty}\dfrac{z^n}{(a+n)^s}$.

It appears that Oberhettinger [45] had some reservations about the rigour employed by Mordell [41].

□

Letting $a=1$ in (3.7.5) results in [36, p.10]

(3.7.11) $$\varsigma(s) = \frac{s}{\pi}\sum_{n=1}^{\infty}\frac{1}{n}\int_0^{\infty}\frac{\sin 2\pi nu}{u^{s+1}}du$$

and with $a=1/2$ in (3.7.5) we obtain

(3.7.12) $$\varsigma\left(s,\frac{1}{2}\right) = \frac{s}{\pi}\sum_{n=1}^{\infty}\frac{(-1)^n}{n}\int_0^{\infty}\frac{\sin 2\pi nu}{u^{s+1}}du$$

Differentiating (3.7.11) gives us

$$\varsigma^{(m)}(0) = (-1)^m\frac{m}{\pi}\sum_{n=1}^{\infty}\frac{1}{n}\int_0^{\infty}\frac{\sin 2\pi nu}{u}\log^m u\,du$$

$$\varsigma^{(m)}\left(0,\frac{1}{2}\right)=(-1)^m \frac{m}{\pi}\sum_{n=1}^{\infty}\frac{(-1)^n}{n}\int_0^{\infty}\frac{\sin 2\pi nu}{u}\log^m u\, du$$

□

Formal derivations of (3.7.6) and (3.7.7) are set out below. Using the definition of the gamma function we have

$$\int_0^{\infty} e^{-zy} y^{s-1} dy = \frac{\Gamma(s)}{z^s}$$

and with $z = u \pm ix$ we obtain

$$\int_0^{\infty} e^{-uy} y^{s-1}[\cos(xy) - i\sin(xy)]dy = (u+ix)^{-s}\Gamma(s)$$

$$\int_0^{\infty} e^{-uy} y^{s-1}[\cos(xy) + i\sin(xy)]dy = (u-ix)^{-s}\Gamma(s)$$

This gives us

$$i\int_0^{\infty} e^{-uy} y^{s-1} \sin(xy) dy = -\left[(u+ix)^{-s} - (u-ix)^{-s}\right]\Gamma(s)$$

$$\int_0^{\infty} e^{-uy} y^{s-1} \cos(xy) dy = \left[(u+ix)^{-s} + (u-ix)^{-s}\right]\Gamma(s)$$

Then, noting that

$$(u+ix)^{-s} - (u-ix)^{-s} = (re^{i\theta})^{-s} - (re^{-i\theta})^{-s}$$

$$= r^{-s}[e^{-is\theta} - e^{is\theta}]$$

we have

$$(u+ix)^{-s} - (u-ix)^{-s} = \frac{2}{i(u^2+x^2)^{s/2}}\sin(s\tan^{-1}(x/u))$$

and we obtain

$$\int_0^{\infty} e^{-uy} y^{s-1} \sin(xy) dy = \Gamma(s)\frac{\sin[s\tan^{-1}(x/u)]}{(u^2+x^2)^{s/2}}$$

and

$$\int_0^\infty e^{-uy} y^{s-1} \cos(xy)\,dy = \Gamma(s)\frac{\cos[s\tan^{-1}(x/u)]}{(u^2+x^2)^{s/2}}$$

Letting $u = 0$ gives us

$$\int_0^\infty y^{s-1} \sin(xy)\,dy = \frac{\Gamma(s)\sin(\pi s/2)}{x^s}$$

and

$$\int_0^\infty y^{s-1} \cos(xy)\,dy = \frac{\Gamma(s)\cos(\pi s/2)}{x^s}$$

which correspond with (3.7.6) and (3.7.7) above.

$\square$

We have the gamma function

$$\Gamma(p) = \int_0^\infty t^{p-1} e^{-t}\,dt$$

and with the substitution $t = ux$ this becomes

$$\frac{1}{x^p} = \frac{1}{\Gamma(p)}\int_0^\infty u^{p-1} e^{-xu}\,du$$

Therefore we obtain upon integration

$$\int_0^\infty \frac{\cos 2\pi nx}{x^p}\,dx = \frac{1}{\Gamma(p)}\int_0^\infty\int_0^\infty u^{p-1} e^{-xu} \cos 2\pi nx\,du\,dx$$

$$= \frac{1}{\Gamma(p)}\int_0^\infty u^{p-1}\,du \int_0^\infty e^{-xu} \cos 2\pi nx\,dx$$

and using

$$\int_0^\infty e^{-ux} \cos 2\pi nx\,dx = \frac{u}{u^2+4\pi^2 n^2}$$

gives us

$$\int_0^\infty \frac{\cos 2\pi nx}{x^p}\,dx = \frac{1}{\Gamma(p)}\int_0^\infty \frac{u^p}{u^2+4\pi^2 n^2}\,du$$

We have the summation

(3.7.13) $$\sum_{n=1}^{\infty}\int_0^{\infty}\frac{\cos 2\pi nx}{x^p}dx = \frac{1}{\Gamma(p)}\sum_{n=1}^{\infty}\int_0^{\infty}\frac{u^p}{u^2+4\pi^2 n^2}du$$

$$= \frac{1}{\Gamma(p)}\int_0^{\infty}\sum_{n=1}^{\infty}\frac{u^p}{u^2+4\pi^2 n^2}du$$

assuming that the interchange of integration and summation is valid.

We then use the identity (3.6.9)

$$2\sum_{n=1}^{\infty}\frac{y}{y^2+(2n\pi)^2} = \frac{1}{e^y-1} - \frac{1}{y} + \frac{1}{2}$$

to obtain

$$\sum_{n=1}^{\infty}\int_0^{\infty}\frac{\cos 2\pi nx}{x^p}dx = \frac{1}{2\Gamma(p)}\int_0^{\infty}\left[\frac{1}{e^u-1}-\frac{1}{u}+\frac{1}{2}\right]u^{p-1}du$$

For $-1 < \sigma < 0$ we have [36, p.24]

$$\varsigma(s) = \frac{1}{\Gamma(p)}\int_0^{\infty}\left[\frac{1}{e^u-1}-\frac{1}{u}+\frac{1}{2}\right]u^{p-1}du$$

and we therefore have

(3.7.14) $$\sum_{n=1}^{\infty}\int_0^{\infty}\frac{\cos 2\pi nx}{x^s}dx = \frac{1}{2}\varsigma(s)$$

From (3.7.4) we see that

$$\varsigma(s,1) = 2\sum_{n=1}^{\infty}\int_{-1}^{\infty}\frac{\cos 2\pi nx}{(1+x)^s}dx$$

$$= 2\sum_{n=1}^{\infty}\int_0^{\infty}\frac{\cos 2\pi n(u-1)}{u^s}du$$

which corresponds with (3.7.14).

□

Integration by parts results in

$$\sum_{n=1}^{\infty} \int_0^{\infty} \frac{\cos 2\pi nx}{x^s} dx = \lim_{N \to \infty} \sum_{n=1}^{\infty} \frac{\sin 2\pi nx}{2\pi n x^s} \bigg|_0^N + s \sum_{n=1}^{\infty} \int_0^{\infty} \frac{\sin 2\pi nx}{2\pi n x^{s+1}} dx$$

$$= s \sum_{n=1}^{\infty} \int_0^{\infty} \frac{\sin 2\pi nx}{2\pi n x^{s+1}} dx$$

and Titchmarsh [50, p.15] tells us that this is equal to $\frac{1}{2}\varsigma(s)$

□

Alternatively, we may consider the integral

$$I = \int_0^{\infty} \frac{u^p}{u^2 + t^2} du = \frac{1}{t^2} \int_0^{\infty} \frac{u^p}{(u/t)^2 + 1} du$$

Letting $v = (u/t)^2$ we get

$$I = \frac{t^{p-1}}{2} \int_0^{\infty} \frac{v^{(p-1)/2}}{1+v} dv$$

$$= \frac{1}{2t^{1-p}} \frac{\pi}{\sin[\pi(p+1)/2]}$$

$$= \frac{1}{2t^{1-p}} \frac{\pi}{\cos(\pi p/2)}$$

where we have used the well known integral

$$\int_0^{\infty} \frac{v^{s-1}}{1+v} dv = \frac{\pi}{\sin \pi s}.$$

Substituting this in (3.7.13) results in (cf [50, p.24])

$$\sum_{n=1}^{\infty} \int_0^{\infty} \frac{\cos 2\pi nx}{x^p} dx = \frac{1}{\Gamma(p)} \sum_{n=1}^{\infty} \frac{1}{2(2\pi n)^{1-p}} \frac{\pi}{\cos(\pi p/2)}$$

$$= \frac{\pi}{2\Gamma(p)(2\pi)^{1-p} \cos(\pi p/2)} \varsigma(1-p)$$

$$=\frac{1}{2}\varsigma(p)$$

where we have employed the functional equation for the Riemann zeta function. This gives us (3.7.14).

□

We have from (3.7.6)

$$\int_0^\infty \frac{\cos u}{u^s}\,du = \frac{\pi}{2\Gamma(s)\cos(\pi s/2)}$$

and the substitution $u = 2\pi nx$ gives us

$$\int_0^\infty \frac{\cos 2\pi nx}{x^s}\,dx = \frac{\pi}{2(2\pi)^{1-s}n^{1-s}\Gamma(s)\cos(\pi s/2)}$$

and we have the summation

$$\sum_{n=1}^\infty \int_0^\infty \frac{\cos 2\pi nx}{x^s}\,dx = \frac{\pi\varsigma(1-s)}{2(2\pi)^{1-s}\Gamma(s)\cos(\pi s/2)}$$

Using the functional equation for the Riemann zeta function

$$\varsigma(1-s) = 2(2\pi)^{-s}\Gamma(s)\cos(\pi s/2)\varsigma(s)$$

we obtain (3.7.14) again

$$\sum_{n=1}^\infty \int_0^\infty \frac{\cos 2\pi nx}{x^s}\,dx = \frac{1}{2}\varsigma(s)$$

□

We see that

$$\int_0^\infty \frac{\cos 2\pi nx}{(1+x)^s}\,dx = \int_1^\infty \frac{\cos 2\pi n(u-1)}{u^s}\,du$$

$$= \int_1^\infty \frac{\cos 2\pi nu}{u^s}\,du$$

and we therefore have

$$\sum_{n=1}^{\infty} \int_0^{\infty} \frac{\cos 2\pi nx}{(1+x)^s} dx = \sum_{n=1}^{\infty} \int_0^{\infty} \frac{\cos 2\pi nx}{x^s} dx - \sum_{n=1}^{\infty} \int_0^1 \frac{\cos 2\pi nx}{x^s} dx$$

$$= \frac{1}{2}\varsigma(s) - \sum_{n=1}^{\infty} \int_0^1 \frac{\cos 2\pi nx}{x^s} dx$$

Letting $a = 1$ in (3.7.1) gives us

$$\varsigma(s) = \frac{1}{2} + \frac{1}{s-1} + 2\sum_{n=1}^{\infty} \int_0^{\infty} \frac{\cos 2\pi nx}{(1+x)^s} dx$$

which appears to suggest that

(3.7.15) $$2\sum_{n=1}^{\infty} \int_0^1 \frac{\cos 2\pi nx}{x^s} dx = \frac{1}{2} + \frac{1}{s-1} + \varsigma(s)$$

which does <u>not</u> agree with equation (8.16) in Berndt's book [10, Part II, p.317]

$$\varsigma(1/2) = \lim_{\varepsilon \to 0+} 2\sum_{n=1}^{\infty} \int_\varepsilon^1 \frac{\cos 2\pi nx}{\sqrt{x}} dx$$

Differentiating (3.7.15) gives us

$$-2\sum_{n=1}^{\infty} \int_0^1 \frac{\log x \cos 2\pi nx}{x^s} dx = -\frac{1}{(s-1)^2} + \varsigma'(s)$$

and with $s = 0$ we have

$$\varsigma'(0) - 1 = -2\sum_{k=1}^{\infty} \int_0^1 \log x \cdot \cos 2\pi kx\, dx$$

It is easily seen that

$$\int_0^1 \log x \cdot \cos 2\pi nx\, dx = \log x \cdot \frac{\sin 2\pi nx}{2\pi n}\Big|_0^1 - \int_0^1 \frac{\sin 2\pi nx}{2\pi nx} dx$$

$$= -\int_0^1 \frac{\sin 2\pi kx}{2\pi kx} dx$$

$$\int_0^1 \log x \cdot \cos 2\pi nx\, dx = -\frac{Si(2\pi n)}{2\pi n}$$

This therefore suggests that

$$\varsigma'(0) - 1 = \sum_{n=1}^{\infty} \frac{Si(2\pi n)}{\pi n}$$

but this does not agree with (3.2.3)

$$\sum_{n=1}^{\infty} \frac{si(2n\pi)}{n} = \frac{\pi}{2}\log(2\pi) - \pi$$

□

We write (3.7.1) as

(3.7.16) $$\varsigma(s,a) - \frac{1}{s-1} = \frac{1}{2}\frac{1}{a^s} + \frac{a^{1-s}-1}{s-1} + 2\sum_{n=1}^{\infty}\int_0^{\infty} \frac{\cos 2\pi nx}{(a+x)^s}dx$$

It may be seen from (3.5.1) that

$$\gamma_p(a) = (-1)^p \lim_{s\to 1} \frac{d^p}{ds^p}\left[\varsigma(s,a) - \frac{1}{s-1}\right]$$

and differentiating (3.7.16) $p$ times gives us

$$\frac{d^p}{ds^p}\left[\varsigma(s,a) - \frac{1}{s-1}\right] = \frac{a^{-s}(-1)^p \log^p a}{2} + f^{(p)}(s) + 2(-1)^p \sum_{n=1}^{\infty}\int_0^{\infty} \frac{\log^p(a+x)\cos 2\pi nx}{(a+x)^s}dx$$

where we have denoted $f(s)$ as

$$f(s) = \frac{a^{1-s}-1}{s-1}$$

We can represent $f(s)$ by the integral

$$f(s) = \frac{a^{1-s}-1}{s-1} = -\int_1^a x^{-s}dx$$

so that

$$f^{(p)}(s) = -(-1)^p \int_1^a x^{-s} \log^p x \, dx$$

and thus

$$f^{(p)}(1) = -(-1)^p \int_1^a \frac{\log^p x}{x} dx = (-1)^{p+1} \frac{\log^{p+1} a}{p+1}$$

Therefore, upon taking the limit $s \to 1$, we obtain for $p \geq 1$

(3.7.17) $$\gamma_p(a) = \frac{1}{2} \frac{\log^p a}{a} - \frac{\log^{p+1} a}{p+1} + 2 \sum_{n=1}^{\infty} \int_0^{\infty} \frac{\log^p(a+x)}{a+x} \cos 2\pi nx \, dx$$

as previously shown in an equivalent form by Zhang and Williams [54, Eq(6.1)] in 1994.

This may be written as

$$\gamma_p(a) = \frac{1}{2} \frac{\log^p a}{a} - \frac{\log^{p+1} a}{p+1} + 2 \sum_{n=1}^{\infty} \int_a^{\infty} \frac{\log^p u}{u} \cos 2\pi n(u-a) \, du$$

and with $a = 1$ we obtain (3.5.5).

Letting $s \to 1$ in (3.7.16) gives us

$$\gamma_0(a) = \frac{1}{2a} - \log a + 2 \sum_{n=1}^{\infty} \int_0^{\infty} \frac{\cos 2\pi nx}{a+x} dx$$

so that

$$\gamma_0(a) = \frac{1}{2a} - \log a - 2 \sum_{n=1}^{\infty} \left[ \cos(2n\pi a) si(2n\pi a) + \sin(2n\pi a) Ci(2n\pi a) \right]$$

and, since $\gamma_0(a) = -\psi(a)$, this is equivalent to (3.7) as reported by Nörlund [44].

□

Differentiation of (3.7.16) with respect to $s$ gives us

$$\varsigma'(s,a) = -\frac{1}{2} \frac{\log a}{a^s} - \frac{(s-1)\log a + 1}{(s-1)^2 a^{s-1}} - 2 \sum_{n=1}^{\infty} \int_0^{\infty} \frac{\cos 2\pi nx \log(a+x)}{(a+x)^s} dx$$

so that

(3.7.18) $$\varsigma'(0,a) = \left(a - \frac{1}{2}\right) \log a - a - 2 \sum_{n=1}^{\infty} \int_0^{\infty} \cos 2\pi nx \log(a+x) \, dx$$

Integration by parts gives us

$$\int_0^N \cos 2\pi nx \log(a+x)\,dx = \frac{\sin 2\pi nx}{2\pi n}\log(a+x)\Big|_0^N - \frac{1}{2\pi n}\int_0^N \frac{\sin 2\pi nx}{a+x}\,dx$$

$$= -\frac{1}{2\pi n}\int_0^N \frac{\sin 2\pi nx}{a+x}\,dx$$

and we therefore have

$$\int_0^\infty \cos 2\pi nx \log(a+x)\,dx = -\frac{1}{2\pi n}\int_0^\infty \frac{\sin 2\pi nx}{a+x}\,dx$$

Using Hurwitz's formula

$$\varsigma'(0,a) = \log\Gamma(a) - \frac{1}{2}\log(2\pi)$$

we obtain

(3.7.19) $\qquad \log\Gamma(a) = \frac{1}{2}\log(2\pi) + \left(a - \frac{1}{2}\right)\log a - a + \frac{1}{\pi}\sum_{n=1}^\infty \int_0^\infty \frac{\sin(2n\pi x)}{n(x+a)}\,dx$

which is one of the exercises posed by Whittaker & Watson [51, p.261]. This result was attributed by Stieltjes to Bourguet. Equation (3.7.19) may also be derived using the Euler-Maclaurin summation formula (see for example Knopp's book [38, p.530]).

With regard to (3.7.18) we have

$$\int_0^N \cos 2\pi nx \log(a+x)\,dx$$

$$= \frac{1}{2n\pi}\Big[\sin 2n\pi a\, Ci\big[2n\pi(a+x)\big] - \cos 2n\pi a\, Si\big[2n\pi(a+x)\big] + \sin 2n\pi x\log(a+x)\Big]\Big|_0^N$$

$$= \frac{1}{2n\pi}\Big[\sin 2n\pi a\, Ci\big[2n\pi(a+x)\big] - \cos 2n\pi a\, Si\big[2n\pi(a+x)\big]\Big]\Big|_0^N$$

Since $\lim_{x\to\infty} Si(x) = \int_0^\infty \frac{\sin t}{t}\,dt = \frac{\pi}{2}$ and $\lim_{x\to\infty} Ci(x) = 0$ we determine that

$$\int_0^\infty \cos 2\pi nx \log(a+x)\,dx = \frac{1}{2n\pi}\left[-\sin 2n\pi a\, Ci(2n\pi a) - \cos 2n\pi a\left\{\frac{\pi}{2} - Si(2n\pi a)\right\}\right]$$

$$= \frac{1}{2n\pi}\left[\cos(2n\pi a)si(2n\pi a) - \sin(2n\pi a)Ci(2n\pi a)\right]$$

since $si(x) = Si(x) - \dfrac{\pi}{2}$.

Therefore we have

$$\varsigma'(0,a) = \left(a - \frac{1}{2}\right)\log a - a - \sum_{n=1}^\infty \frac{1}{n\pi}\left[\cos 2n\pi a\, si(2n\pi a) - \sin 2n\pi a\, Ci(2n\pi a)\right]$$

which implies that

$$\log \Gamma(a) =$$

$$\frac{1}{2}\log(2\pi) + \left(a - \frac{1}{2}\right)\log a - a + \frac{1}{\pi}\sum_{n=1}^\infty \frac{1}{n}[\sin(2n\pi a)Ci(2n\pi a) - \cos(2n\pi a)si(2n\pi a)]$$

which we saw earlier in (3.10).

We have the second derivative of (3.7.16)

$$\varsigma''(s,a) = \frac{1}{2}\frac{\log^2 a}{a^s} + \frac{2(s-1)\log a + (s-1)^2\log^2 a + 2}{(s-1)^3 a^{s-1}} + 2\sum_{n=1}^\infty \int_0^\infty \frac{\cos 2\pi nx \log^2(a+x)}{(a+x)^s}\,dx$$

so that

$$\varsigma''(0,a) = -\left(a - \frac{1}{2}\right)\log^2 a + 2a\log a - 2 + 2\sum_{n=1}^\infty \int_0^\infty \cos 2\pi nx \log^2(a+x)\,dx$$

Integration by parts gives us

$$\int_0^N \cos 2\pi nx \log^2(a+x)\,dx = \frac{\sin 2\pi nx}{2\pi n}\log^2(a+x)\bigg|_0^N - \frac{1}{\pi n}\int_0^N \frac{\sin 2\pi nx \log(a+x)}{a+x}\,dx$$

$$= -\frac{1}{\pi n}\int_0^N \frac{\sin 2\pi nx \log(a+x)}{a+x}\,dx$$

$$-\frac{1}{\pi n}\int_0^N \frac{\sin 2\pi nx \log(a+x)}{a+x}dx = \frac{\cos 2\pi nx \log(a+x)}{2(\pi n)^2(a+x)}\bigg|_0^N + \frac{1}{2(\pi n)^2}\int_0^N \frac{\cos 2\pi nx[1-\log(a+x)]}{(a+x)^2}dx$$

$$= \frac{\log(a+N)}{2(\pi n)^2(a+N)} - \frac{\log a}{2(\pi n)^2 a} + \frac{1}{2(\pi n)^2}\int_0^N \frac{\cos 2\pi nx[1-\log(a+x)]}{(a+x)^2}dx$$

and, since $\lim_{N\to\infty}\frac{\log(a+N)}{(a+N)} = 0$, we have

$$\varsigma''(0,a) = -\left(a-\frac{1}{2}\right)\log^2 a + 2a\log a - 2 - \frac{\varsigma(2)\log a}{\pi^2 a} + \frac{1}{\pi^2}\sum_{n=1}^\infty \int_0^\infty \frac{\cos 2\pi nx[1-\log(a+x)]}{n^2(a+x)^2}dx$$

We note that Mathematica evaluates $\int \cos 2\pi nx \log^2(1+x)\,dx$ in terms of the generalised hypergeometric function $_3F_3$ and the incomplete gamma function.

### 4. Some integrals involving $\cot x$

In this section we shall have need of Bonnet's second mean-value theorem which states [17, p.110]:

(i) Let $\phi(x)$ be bounded, monotonic decreasing and never negative in $[a,b]$; and let $\psi(x)$ be bounded and integrable in $[a,b]$. Then we have

$$\int_a^b \phi(x)\psi(x)dx = \phi(a)\int_a^\xi \psi(x)dx$$

where $\xi$ is some definite value of $x$ in $a \le \xi \le b$.

(ii) Alternatively, subject to the same conditions as above but with $\phi(x)$ monotonic increasing, we then have

$$\int_a^b \phi(x)\psi(x)dx = \phi(b)\int_\xi^b \psi(x)dx$$

We now use Bonnet's second mean-value theorem to prove the following proposition.

**Proposition**

We assume that $f(x)$ is continuous on $[0,b]$ and that $f(0) = 0$. Let us also assume that $f(x)$ is non-negative and monotonic increasing on $[0,b]$.

Then we have

(4.1) $$\int_0^b f(x)\cot(x/2)\,dx = 2\sum_{n=1}^{\infty} \int_0^b f(x)\sin 2nx\,dx$$

**Proof**

We recall (1.3)

$$\frac{1}{2}\cot(x/2) = \sum_{n=0}^{N} \sin nx + \frac{\cos(N+1/2)x}{2\sin(x/2)}$$

$$= \sum_{n=1}^{N} \sin nx + \cot(x/2)\cos Nx - \sin Nx$$

We now multiply this by $f(x)$ and integrate over the interval $[\pi/N, b]$ where $b \leq \pi$.

Provided the relevant integrals exist, we easily see that

$$\frac{1}{2}\int_{\pi/N}^{b} f(x)\cot(x/2)\,dx = \sum_{n=1}^{N} \int_{\pi/N}^{b} f(x)\sin nx\,dx + \int_{\pi/N}^{b} f(x)\cot(x/2)\cos Nx\,dx - \int_{\pi/N}^{b} f(x)\sin Nx\,dx$$

$$= \sum_{n=1}^{N} J_n + I_1(N) + I_2(N)$$

and, first of all, we employ the approach adopted by Zygmund [55, p.59] for the integral

$$I_1(N) = \int_{\pi/N}^{b} f(x)\cot(x/2)\cos Nx\,dx.$$

Since we assumed that $f(x)$ is continuous and that $f(0) = 0$, given $\varepsilon > 0$, we may choose $\eta > 0$ such that $f(\eta) < \varepsilon$ and we write

$$\int_{\pi/N}^{b} f(x)\cot(x/2)\cos Nx\,dx = \int_{\pi/N}^{\eta} f(x)\cot(x/2)\cos Nx\,dx + \int_{\eta}^{b} f(x)\cot(x/2)\cos Nx\,dx$$

and we have the inequality

$$\left|\int_{\pi/N}^{b} f(x)\cot(x/2)\cos Nx\,dx\right| \le \left|\int_{\pi/N}^{\eta} f(x)\cot(x/2)\cos Nx\,dx\right| + \left|\int_{\eta}^{b} f(x)\cot(x/2)\cos Nx\,dx\right|$$

We note that $\cot(x/2)$ is monotonic decreasing on $[\pi/N,\pi]$ and never negative in $[\pi/N,\pi]$ (which is the reason why we specified that $b \le \pi$) and hence the second mean-value theorem (i) gives us for the first part

(4.2) $$\int_{\pi/N}^{\eta} \cot(x/2)f(x)\cos Nx\,dx = \cot(\pi/2N)\int_{\pi/N}^{\eta'} f(x)\cos Nx\,dx$$

Since we assumed that $f(x)$ is non-negative and monotonic increasing on $[0,b]$, we may apply the second limb of Bonnet's second mean-value theorem (ii) to the integral on the right-hand side of (4.2)

$$\int_{\pi/N}^{\eta'} f(x)\cos Nx\,dx = f(\eta')\int_{\eta''}^{\eta'} \cos Nx\,dx$$

and we obtain

$$\int_{\pi/N}^{\eta} f(x)\cot(x/2)\cos Nx\,dx = \cot(\pi/2N)f(\eta')\int_{\eta''}^{\eta'} \cos Nx\,dx$$

$$= \cot(\pi/2N)f(\eta')\frac{[\sin N\eta' - \sin N\eta'']}{N}$$

We recall the elementary inequality $\sin x < x < \tan x$ for $0 < x < \pi/2$ which yields

$$0 < \cot x < \frac{1}{x}$$

and we then have

$$\left|\int_{\pi/N}^{\eta} f(x)\cot(x/2)\cos Nx\,dx\right| < \frac{2N}{\pi}\varepsilon\frac{2}{N}$$

Since $\eta > 0$, the strong version of the Riemann-Lebesgue lemma [5, p.313] tells us that

$$\lim_{N\to\infty} \int_{\eta}^{b} f(x)\cot(x/2)\cos Nx\, dx = 0$$

and similarly we also have $\lim_{N\to\infty} I_2(N) = 0$. Therefore as $N \to \infty$ we have the limit

$$\int_0^b f(x)\cot(x/2)\,dx = 2\sum_{n=1}^{\infty} \int_0^b f(x)\sin 2nx\,dx$$

where, as stated above, we require that $f(0) = 0$. We may also show that

(4.3) $$\int_0^b f(x)\cot(\alpha x/2)\,dx = 2\sum_{n=1}^{\infty} \int_0^b f(x)\sin 2n\alpha x\,dx$$

It is easily seen that this may be generalised to

(4.4) $$\int_a^b f(x)\cot(\alpha x/2)\,dx = 2\sum_{n=1}^{\infty} \int_a^b f(x)\sin 2n\alpha x\,dx$$

We therefore have a version of (1.4a) which now encompasses a larger class of eligible functions $f(x)$.

$\square$

Various applications of (1.4a) were considered in [] and these included the evaluation of the following integrals:

$$\int_0^{\pi/6} x^2 \cot x\, dx = -\frac{1}{3}\varsigma(3) + \frac{\pi^2}{36}\log[2\sin(\pi/6)] + \frac{\sqrt{3}}{6}\pi\left[-\frac{4}{9}\varsigma(2) + \frac{1}{36}\left\{\varsigma\left(2,\frac{1}{6}\right) + \varsigma\left(2,\frac{1}{3}\right)\right\}\right]$$

$$\int_0^{\pi/8} x\cot x\, dx = \frac{\pi}{16}\log\left[2-\sqrt{2}\right] + \frac{1}{8}\left[1-\sqrt{2}\right]G + \frac{1}{64}\left[\sqrt{2}\varsigma\left(2,\frac{1}{8}\right) - 2\left(\sqrt{2}+1\right)\pi^2\right]$$

$$\frac{1}{2}\int_0^1 (x^2\log x)\cot(\pi x/2)\,dx =$$

$$-\frac{3}{\pi^3}[\varsigma_a(3) + \varsigma(3)] - \frac{2}{\pi^3}\sum_{n=1}^{\infty}\frac{Ci(n\pi)}{n^3} + \frac{2}{\pi^3}\varsigma(3)[\gamma + \log\pi] - \frac{2}{\pi^3}\varsigma'(3)$$

$$\int_0^1 \log \Gamma(x+1)\cot \pi x\, dx = 2\sum_{n=1}^{\infty}\int_0^1 \log \Gamma(x+1)\sin 2n\pi x\, dx = \frac{1}{\pi}\sum_{n=1}^{\infty}\frac{Ci(2n\pi)}{n}$$

□

Since $\psi(x)-\psi(1-x)=-\pi\cot \pi x$ we see that

(4.5) $$\int_0^1 f(x)[\psi(1-x)-\psi(x)]\,dx = 2\pi\sum_{n=1}^{\infty}\int_0^1 f(x)\sin 2\pi nx\, dx$$

Part of the following is based on an observation made by Glasser [29] in 1966. Let us consider the integral

$$I = \int_0^1 f(x)\psi(x)\,dx$$

where $f(x) = -f(1-x)$ and $f(x)$ is selected so that the integral converges. We then have

$$I = \int_0^1 f(1-t)\psi(1-t)\,dt = -\int_0^1 f(t)\psi(1-t)\,dt$$

and hence we see that

$$2I = \int_0^1 f(x)[\psi(x)-\psi(1-x)]\,dx$$

Therefore, since $\psi(x)-\psi(1-x)=-\pi\cot \pi x$, we have

(4.6) $$\int_0^1 f(x)\psi(x)\,dx = -\frac{\pi}{2}\int_0^1 f(x)\cot \pi x\, dx = -\pi \int_0^{1/2} f(x)\cot \pi x\, dx$$

Integration by parts formally gives us

$$I = \int_0^1 f(x)\psi(x)\,dx = f(x)\log \Gamma(x)\Big|_0^1 - \int_0^1 f'(x)\log \Gamma(x)\,dx$$

and, for suitably behaved functions, we have

$$\int_0^1 f(x)\psi(x)\,dx = -\int_0^1 f'(x)\log\Gamma(x)\,dx$$

Hence we see that

(4.7) $\quad \dfrac{\pi}{2}\int_0^1 f(x)\cot\pi x\,dx = \pi\int_0^{1/2} f(x)\cot\pi x\,dx = \int_0^1 f'(x)\log\Gamma(x)\,dx$

where $f(x) = -f(1-x)$ and $f(x)$ is selected so that the integral converges.

An example of this is

$$\int_0^1 B_{2n+1}(x)\psi(x)\,dx = B_{2n+1}(x)\log\Gamma(x)\Big|_0^1 - (2n+1)\int_0^1 B_{2n}(x)\log\Gamma(x)\,dx$$

$$= -(2n+1)\int_0^1 B_{2n}(x)\log\Gamma(x)\,dx$$

Glasser's formula (4.6) gives us

(4.8) $\quad \displaystyle\int_0^1 B_{2n+1}(x)\psi(x)\,dx = -\dfrac{\pi}{2}\int_0^1 B_{2n+1}(x)\cot\pi x\,dx$

and hence we have

(4.9) $\quad \displaystyle\int_0^1 B_{2n}(x)\log\Gamma(x)\,dx = \dfrac{\pi}{2(2n+1)}\int_0^1 B_{2n+1}(x)\cot\pi x\,dx$

The following integral appears in Abramowitz and Stegun [1, p.807]

(4.10) $\quad \displaystyle\int_0^1 B_{2n+1}(x)\cot\pi x\,dx = \dfrac{(-1)^{n+1}2(2n+1)!\varsigma(2n+1)}{(2\pi)^{2n+1}}$

and there are many derivations of this; for example, see the recent one by Dwilewicz and Mináč [27]. A similar identity was also derived by Espinosa and Moll [28] in the form

(4.11) $\quad \displaystyle\int_0^1 B_{2n}(x)\log\sin\pi x\,dx = (-1)^n\dfrac{(2n)!\varsigma(2n+1)}{(2\pi)^{2n}}$

and it is easily shown that equation (4.11) above is equivalent to (4.10) following a simple integration by parts.

Hence we obtain

$$(4.12) \quad \int_0^1 B_{2n}(x) \log \Gamma(x)\, dx = \frac{(-1)^{n+1}(2n)!\varsigma(2n+1)}{2(2\pi)^{2n+1}}$$

in agreement with [28].

**5. Open access to our own work**

This paper contains references to various other papers and, rather surprisingly, most of them are currently freely available on the internet. Surely now is the time that <u>all</u> of <u>our</u> work should be freely accessible by <u>all</u>. The mathematics community should lead the way on this by publishing <u>everything</u> on arXiv, or in an equivalent open access repository. We think it, we write it, so why hide it? You know it makes sense.

**REFERENCES**

[1]   M. Abramowitz and I.A. Stegun (Eds.), Handbook of Mathematical Functions with Formulas, Graphs and Mathematical Tables. Dover, New York, 1970.
      http://www.math.sfu.ca/~cbm/aands/

[2]   V.S. Adamchik and H.M. Srivastava, Some Series of the Zeta and Related Functions. Analysis 18, 131-144, 1998.
      http://www-2.cs.cmu.edu/~adamchik/articles/sums.htm

[3]   V.S. Adamchik, Symbolic and numeric computations of the Barnes function.
      Computer Physics Comms., 157 (2004) 181-190.
      http://www.cs.cmu.edu/~adamchik/articles/PhysCom.pdf

[4]   V.S. Adamchik, On the Hurwitz function for rational arguments.
      Applied Mathematics and Computation, 187 (2007) 3-12.
      http://www.cs.cmu.edu/~adamchik/articles/AMC.pdf

[5]   T.M. Apostol, Mathematical Analysis, Second Ed., Addison-Wesley Publishing Company, Menlo Park (California), London and Don Mills (Ontario), 1974.

[6]   H. Azuma, Application of Abel-Plana formula for collapse and revival of Rabi oscillations in Jaynes-Cummings model.
      Int. J. Mod. Phys. C 21, 1021-1049 (2010) arXiv:0901.4357 [pdf, ps, other]

[7]   G. Bachman, L. Narici and E. Beckenstein, Fourier and Wavelet Analysis.
      Springer-Verlag, New York, 2000.

[8]   N. Batir, Inequalities for the double gamma function. 2012.
      http://ajmaa.org/RGMIA/papers/v11n4/JMAA_BARNES.pdf

[9]    B.C. Berndt, On the Hurwitz zeta function.
      Rocky Mtn. J. Math. 2, 151-157, (1972).

[10]  B.C. Berndt, Ramanujan's Notebooks. Part II, Springer-Verlag, 1989.

[11]  B.C. Berndt, Ramanujan's Notebooks. Part V, Springer-Verlag, 1998.

[12]  B.C. Berndt and A. Dixit, A transformation formula involving the Gamma and
      Riemann zeta functions in Ramanujan's Lost Notebook.
      arXiv:0904.1053 [pdf, ps, other], 2009.

[13]  M.T. Boudjelkha, A proof that extends Hurwitz formula into the critical strip.
      Applied Mathematics Letters, 14 (2001) 309-403.

[14]  W.E. Briggs, Some constants associated with the Riemann zeta-function.
      (1955-1956), Michigan Math. J. 3, 117-121.

[15]  T.J.I'a Bromwich, Introduction to the theory of infinite series. Third edition.
      AMS Chelsea Publishing, 1991.

[16]  P. L. Butzer, P. J. S. G. Ferreira, G. Schmeisser and R. L. Stens,
      The Summation Formulae of Euler–Maclaurin, Abel–Plana, Poisson, and their
      interconnections with the Approximate Sampling Formula of Signal Analysis.
      Results. Math., 59, No.3-4, 359-400, 2011, DOI 10.1007/s00025-010-0083-8
      http://www.ieeta.pt/~pjf/PDF/Ferreira2011b.pdf

[17]  H.S. Carslaw, Introduction to the theory of Fourier Series and Integrals.
      Third Ed. Dover Publications Inc, 1930.

[18]  H. Cohen, Number Theory. Volume II: Analytic and modern tools.
      Springer Science, 2007.

[19]  D.F. Connon, Some series and integrals involving the Riemann zeta function,
      binomial coefficients and the harmonic numbers. Volume I, 2007.
      arXiv:0710.4022 [pdf]

[20]  D.F. Connon, Some series and integrals involving the Riemann zeta function,
      binomial coefficients and the harmonic numbers. Volume V, 2007.
      arXiv:0710.4047 [pdf]

[21]  D.F. Connon, Some series and integrals involving the Riemann zeta function,
      binomial coefficients and the harmonic numbers. Volume VI, 2007.
      arXiv:0710.4032 [pdf]

[22]   D.F. Connon, Some applications of the Stieltjes constants.
       arXiv:0901.2083 [pdf], 2009.

[23]   D.F. Connon, Some trigonometric integrals involving $\log \Gamma(x)$ and the digamma
       function. arXiv:1005.3469 [pdf], 2010.

[24]   D.F. Connon, Some integrals involving the Stieltjes constants: Part II.
       arXiv:1104.1911 [pdf], 2011.

[25]   D.F. Connon, Some applications of the sine and cosine integrals. 2012.

[26]   P.G.L. Dirichlet, Sur la convergence des séries trigonométriques qui servent à
       représenter une fonction arbitraire entre des limites données.
       Journal für die reine und angewandte Mathematik, 4, (1829), 157-169.
       http://www.digizeitschriften.de/zeitschriften/open-access/

[27]   R. Dwilewicz and J. Mináč, An introduction to relations between the
       values of $\varsigma(s)$ in terms of holomorphic functions of two variables.
       Proceedings of the Hayama Symposium on Several Complex Variables,
       Japan, Dec. 2000. Pages 28-38 (2001).
       (see also http://www.mat.uab.es/matmat/PDFv2009/v2009n06.pdf)

[28]   O. Espinosa and V.H. Moll, On some integrals involving the Hurwitz zeta
       function: Part I. The Ramanujan Journal, 6,150-188, 2002.
       http://arxiv.org/abs/math.CA/0012078

[29]   M.L. Glasser, Evaluation of some integrals involving the $\psi$ - function.
       Math. of Comp., Vol.20, No.94, 332-333, 1966.

[30]   I.S. Gradshteyn and I.M. Ryzhik, Tables of Integrals, Series and Products.
       Sixth Ed., Academic Press, 2000.
       Errata for Sixth Edition http://www.mathtable.com/errata/gr6_errata.pdf

[31]   A.P. Guinand, On Poisson's summation formula.
       Annals of Mathematics, Vol. 42, No. 3 (Jul., 1941), pp. 591-603.

[32]   A.P. Guinand, Some formulae for the Riemann zeta-function.
       J. London Math. Soc. 22 (1947), 14–18.

[33]   A.P. Guinand, Some finite identities connected with Poisson's summation formula.
       Proc. Edinburgh Math. Soc. (2) 12 (1960), 17–25.

[34]   G.H. Hardy, Divergent Series. Chelsea Publishing Company, New York, 1991.

[35]   C. Hermite, Correspondance d'Hermite et de Stieltjes. Gauthier-Villars, Paris,

1905. http://ebooks.library.cornell.edu/m/math/browse/author/a.php

[36] A. Ivić, The Riemann Zeta- Function: Theory and Applications.
Dover Publications Inc, 2003.

[37] C. Jordan, Sur la série de Fourier.
Comptes Rendus Acad. Sciences, 92, 228-230, 1881.
http://gallica.bnf.fr/ark:/12148/bpt6k7351t/f227.image.r=camille+jordan.langFR

[38] K. Knopp, Theory and Application of Infinite Series.
Second English Ed. Dover Publications Inc, New York, 1990.

[39] E.E. Kummer, Beitrag zur Theorie der Function $\Gamma(x) = \int_0^\infty e^{-v} v^{x-1} dv$.

J. Reine Angew. Math., 35, 1-4, 1847.
http://www.digizeitschriften.de/

[40] M. Merkle and M. M. R. Merkle, Krull's theory for the double gamma function.
Appl. Math. Comput., 218 (2011), 935-943, doi:10.1016/j.amc.2011.01.090.
http://www.milanmerkle.com/documents/radovi/KRULLcs.pdf

[41] L.J. Mordell, Some applications of Fourier series in the analytic theory
of numbers, Proc. Cambridge Phil. Soc., 24 (1928), 585-596.

[42] L.J. Mordell, Poisson's Summation Formula and the Riemann Zeta Function.
J. London Math. Soc. (1929) s1-4 (4): 285-291.

[43] N. Nielsen, Die Gammafunktion. Chelsea Publishing Company, Bronx and
New York, 1965.

[44] N.E. Nörlund, Vorlesungen über Differenzenrechnung.Chelsea, 1954.
http://math-doc.ujf-grenoble.fr/cgi-bin/linum?aun=001355

[45] F. Oberhettinger, Note on the Lerch zeta function.
Pacific J. Math., Volume 6, Number 1 (1956), 117-120.
http://projecteuclid.org/DPubS?service=UI&version=1.0&verb=Display&handle=euclid.pjm/1103044247

[46] A.A. Saharian, The generalized Abel-Plana formula with applications to Bessel
functions and Casimir effect. arXiv:0708.1187 [pdf, ps, other], 2008.

[47] J. Sondow, Double Integrals for Euler's Constant and $\log \frac{4}{\pi}$ and an Analog
of Hadjicostas's Formula. Amer. Math. Monthly, 112, 61-65, 2005.
arXiv:math/0211148 [pdf]

[48] M.R. Spiegel, Schaum's Outline Series of Theory and Problems of Advanced

Calculus. McGraw Hill Book Company, 1963.

[49]  H.M. Srivastava and J. Choi, Series Associated with the Zeta and Related Functions. Kluwer Academic Publishers, Dordrecht, the Netherlands, 2001.

[50]  E.C. Titchmarsh, The Theory of Functions.
$2^{nd}$ Ed., Oxford University Press, 1932.

[51]  E.T. Whittaker and G.N. Watson, A Course of Modern Analysis: An Introduction to the General Theory of Infinite Processes and of Analytic Functions; With an Account of the Principal Transcendental Functions. Fourth Ed., Cambridge University Press, Cambridge, London and New York, 1963.

[52]  Z.X. Wang and D.R. Guo, Special Functions.
World Scientific Publishing Co Pte Ltd, Singapore, 1989.

[53]  J.R. Wilton, A Proof of Poisson's Summation Formula.
J. London Math. Soc. (1930) s1-5(4): 276-279.

[54]  N.-Y. Zhang and K.S. Williams, Some results on the generalized Stieltjes constants. Analysis 14, 147-162 (1994).
http://people.math.carleton.ca/~williams/papers/pdf/187.pdf

[55]  A. Zygmund, Trigonometric Sums. Cambridge Mathematical Library, 2002.

Wessex House,
Devizes Road,
Upavon,
Pewsey,
Wiltshire SN9 6DL